\documentclass[11pt,twoside]{article}
\newcommand{\documentdate}{20 September 2017}
\usepackage{a4wide,latexsym,graphicx,varioref}

\newcommand{\numsection}[1]{\section{#1}\setcounter{equation}{0}}
\newcommand{\appnumsection}[1]{\section*{#1}\setcounter{equation}{0}
\renewcommand{\theequation}{A.\arabic{equation}}
\renewcommand{\thetheorem}{A.\arabic{theorem}}
\renewcommand{\thetable}{A.\arabic{table}}
\renewcommand{\thefigure}{A.\arabic{figure}}
\renewcommand{\thesection}{A} }

\renewcommand{\theequation}{\arabic{section}.\arabic{equation}}
\renewcommand{\thetable}{\arabic{section}.\arabic{table}}
\renewcommand{\thefigure}{\arabic{section}.\arabic{figure}}

\newcommand{\beqn}[1]{\begin{equation}\label{#1}}
\newcommand{\eeqn}{\end{equation}}

\newcommand{\req}[1]{(\ref{#1})}

\newcommand{\bpr}{{\bf Proof.} \hspace{1.5mm}}
\newcommand{\epr}{\hfill $\Box$ \vspace*{1em}}

\newcommand{\proof}[1]{
\begin{list}{}{
\setlength{\topsep}{0.0pt}
\setlength{\partopsep}{0.0pt}
\setlength{\leftmargin}{0.025\textwidth}
\setlength{\rightmargin}{0.5\leftmargin}
\setlength{\labelwidth}{0.5\leftmargin}
\setlength{\labelsep}{0.25\leftmargin}}
\item \bpr #1 \epr \noindent
\end{list}}

\newcommand{\ms}{\;\;\;\;}
\newcommand{\tim}[1]{\;\; \mbox{#1} \;\;}

\newtheorem{theorem}{Theorem}[section]
\newtheorem{lemma}[theorem]{Lemma}

\newcommand{\llem}[2]{\vspace{\baselineskip} 
\noindent\framebox[\textwidth]{\parbox{0.95\textwidth}{
\begin{lemma} \label{#1} \rm #2 \end{lemma} } } \vspace{\baselineskip} }

\newcommand{\lthm}[2]{\vspace{\baselineskip} 
\noindent\framebox[\textwidth]{\parbox{0.95\textwidth}{
\begin{theorem} \label{#1} \rm #2 \end{theorem} } } \vspace{\baselineskip} }

\newcommand{\calK}{{\cal K}}
\newcommand{\calM}{{\cal M}}
\newcommand{\calO}{{\cal O}}
\newcommand{\calS}{{\cal S}} 

\renewcommand{\Re}{\hbox{I\hskip -2pt R}}

\newcommand{\sfrac}[2]{{\scriptstyle \frac{#1}{#2}}}
\newcommand{\half}{\sfrac{1}{2}}

\newcommand{\quarter}{\sfrac{1}{4}}

\newcommand{\tenth}{\sfrac{1}{10}}

\newcommand{\comment}[1]{}
\newcommand{\eqdef}{\stackrel{\rm def}{=}}
\newcommand{\bigfrac}[2]{\frac{\displaystyle #1}{\displaystyle #2}}
\newcommand{\bigsum}{\displaystyle \sum}
\newcommand{\ii}[1]{\{1, \ldots, #1 \}}
\newcommand{\iibe}[2]{\{ #1, \ldots, #2 \}}

\newcommand{\paperauthor}{Coralia Cartis\hspace*{-0.02cm}\footnotemark[1],\,\, 
Nicholas I.~M.~Gould\footnotemark[2]\,\, and\,
Philippe L.~Toint\footnotemark[3]}
\newcommand{\papertitle}{Worst-case evaluation complexity and optimality of
  second-order methods for nonconvex smooth optimization}

\title{\papertitle}
\author{\paperauthor}
\date{\documentdate}

\begin{document}
\maketitle

\footnotetext[1]{
       Mathematical Institute,
       Oxford University,
       Oxford OX2 6GG, England, United Kingdom.
       Email: coralia.cartis@maths.ox.ac.uk},
\footnotetext[2]{Scientific Computing Department,
       STFC-Rutherford Appleton Laboratory, 
       Chilton, Oxfordshire, OX11 0QX, England, United Kingdom. 
       Email: nick.gould@stfc.ac.uk.
       This work was supported by the EPSRC grant EP/M025179/1.} 
\footnotetext[3]{[Speaker]
       Department of Mathematics,
       University of Namur,
       61, rue de Bruxelles, B-5000, Namur, Belgium.
       Email: philippe.toint@unamur.be.}

\begin{abstract}
We establish or refute the optimality of inexact second-order methods for
unconstrained nonconvex optimization from the point of view of worst-case
evaluation complexity, improving and generalizing the results of
\cite{CartGoulToin10a,CartGoulToin11c}. To this aim, we consider a new general
class of inexact second-order algorithms for unconstrained optimization that
includes regularization and trust-region variations of Newton's method as well
as of their linesearch variants.  For each method in this class and arbitrary
accuracy threshold $\epsilon \in (0,1)$, we exhibit a smooth objective
function with bounded range, whose gradient is globally Lipschitz continuous
and whose Hessian is $\alpha-$H\"older continuous (for given $\alpha\in
[0,1]$), for which the method in question takes at least
$\lfloor\epsilon^{-(2+\alpha)/(1+\alpha)}\rfloor$ function evaluations to
generate a first iterate whose gradient is smaller than $\epsilon$ in norm.
Moreover, we also construct another function on which Newton's takes
$\lfloor\epsilon^{-2}\rfloor$ evaluations, but whose Hessian is Lipschitz
continuous on the path of iterates. These examples provide lower bounds on the
worst-case evaluation complexity of methods in our class when applied to
smooth problems satisfying the relevant assumptions. Furthermore, for
$\alpha=1$, this lower bound is of the same order in $\epsilon$ as the upper
bound on the worst-case evaluation complexity of the cubic
%(and more generally $(2+\alpha)$)
regularization method and other methods in a class of methods
proposed in \cite{CurtRobiSama17c} or in \cite{RoyeWrig17}, thus implying that
these methods have optimal worst-case evaluation complexity within a wider
class of second-order methods, and that Newton's method is suboptimal.
%the same being true for their linesearch variants.
\end{abstract}

\setcounter{page}{1}
\numsection{Introduction}
\label{intro0}

Newton's method has long represented a benchmark for rapid asymptotic
convergence when minimizing smooth, unconstrained objective functions
\cite{DennSchn83}.  It has also been efficiently safeguarded to ensure its
global convergence to first- and even second-order critical points, in the
presence of local nonconvexity of the objective using linesearch
\cite{NoceWrig99}, trust-region \cite{ConnGoulToin00} or other regularization
techniques \cite{Grie81, NestPoly06,CartGoulToin11}.  Many variants of these
globalization techniques have been proposed.  These generally retain fast
local convergence under non-degeneracy assumptions, are often suitable when
solving large-scale problems and sometimes allow approximate rather than true
Hessians to be employed.  We attempt to capture the common features of these
methods in the description of a general class of second-order methods, which
we denote by $\calM.\alpha$ in what follows.

In this paper, we are concerned with establishing {\it lower bounds} on the
worst-case evaluation complexity of the $\calM.\alpha$ methods\footnote{And, as an
aside, on that of the steepest-descent method.} when applied to
``sufficiently smooth'' nonconvex minimization problems, in the sense that we
exhibit objective functions on which these methods take a large number of
function evaluations to obtain an approximate first-order point.

There is a growing literature on the global worst-case evaluation complexity
of first- and second-order methods for nonconvex smooth optimization problems
(for which we provide a partial bibliography with this paper).
\nocite{Agaretal16,AnanGe16,BergDiouGrat17,BianChen13,BianChen15,BianChenYe15,
  BianLiuzMoriScia15,BianScia16,
  BirgGardMartSantToin17,BirgGardMartSantToin16,BirgMart17,BoumAbsiCart16,
  CarmDuch16,CarmDuchHindSidf17,
  CartGoulToin11,CartGoulToin11b,CartGoulToin11c,CartGoulToin11d,CartGoulToin12e,
  CartGoulToin13a,CartGoulToin10a,CartGoulToin12b,CartGoulToin15b,
  CartGoulToin17a,CartGoulToin17b,CartGoulToin17c,CartGoulToin17d,
  CartSampToin15,CartSche17,ChenToinWang17,
  CurtRobiSama17,CurtRobiSama17b,CurtRobiSama17c,
  DodaViceZhan15,Duss15,DussOrba17,FaccKungLampScut17,
  GarmJudiVice16,GhadLan16,GeJianYe11,
  GrapYuanYuan15a,GrapYuanYuan16,GratRoyeViceZhan15,GratRoyeVice17,
  GratSartToin08,GoulPorcToin12,Jarr13,JianLinMaZhan16,Hong06,
  LuWeiLi12,Mart17,MartRayd16,NestPoly06,NestGrap16,
  RoyeWrig17,ScheTang13,ScheTang16,UedaYama10,UedaYama10b,
  Vava93,Vice13,XuRoosMaho17}
In particular, it is known \cite{Vava93}, \cite[p.~29]{Nest04} that
steepest-descent method with either exact or inexact linesearches takes at
most\footnote{When $\{a_k\}$ and $\{b_k\}$ are two sequences of real numbers,
we say that $a_k = \calO\left(b_k\right)$ if the ratio $a_k/b_k$ is bounded.}
$\calO\left(\epsilon^{-2}\right)$ iterations/function-evaluations to generate a
gradient whose norm is at most $\epsilon$ when started from an arbitrary
initial point and applied to nonconvex smooth objectives with gradients that
are globally Lipschitz continuous within some open convex set containing the
iterates generated.  Furthermore, this bound is essentially sharp (for inexact
\cite{CartGoulToin10a} and exact \cite{CartGoulToin12g} linesearches).
Similarly, trust-region methods that ensure at least a Cauchy
(steepest-descent-like) decrease on each iteration satisfy a worst-case
evaluation complexity bound of the same order under identical conditions
\cite{GratSartToin08}.  It follows that Newton's method globalized by
trust-region regularization has the same $\calO\left(\epsilon^{-2}\right)$
worst-case evaluation upper bound; such a bound has also been shown to be
essentially sharp \cite{CartGoulToin10a}.

From a worst-case complexity point of view, one can do better when a cubic
regularization/perturbation of the Newton direction is used \cite{Grie81,
  NestPoly06, CartGoulToin11,CurtRobiSama17c}---such a method iteratively
calculates step corrections by (exactly or approximately) minimizing a cubic
model formed of a quadratic approximation of the objective and the cube of a
weighted norm of the step.  For such a method, the worst-case global
complexity improves to be $\calO\left(\epsilon^{-3/2}\right)$ \cite{NestPoly06,
  CartGoulToin11}, for problems whose gradients and Hessians are Lipschitz
continuous as above; this bound is also essentially sharp
\cite{CartGoulToin10a}.  If instead powers between two and three are used in
the regularization, then an ``intermediate'' worst-case complexity of
$\calO\left(\epsilon^{-(2+\alpha)/(1+\alpha)}\right)$ is obtained for such variants
when applied to functions with globally $\alpha-$H\"older continuous Hessian
on the path of iterates, where $\alpha \in (0,1]$ \cite{CartGoulToin11c}.  It
is finally possible, as proposed in \cite{RoyeWrig17}, to obtain the desired
$\calO\left(\epsilon^{-3/2}\right)$ order of worst-case evaluation complexity
using a purely quadratic regularization, at the price of mixing iterations
using the regularized and unregularized  Hessian with iterations requiring the
computation of its left-most eigenpair.

These (essentially tight) upper bounds on the worst-case evaluation complexity
of such second-order methods naturally raise the question as to whether other
second-order methods might have better worst-case complexity than cubic (or
similar) regularization over certain classes of sufficiently smooth
functions. To attempt to answer this question, we define a general,
parametrized class of methods that includes Newton's method, and that attempts
to capture the essential features of globalized Newton variants we have
mentioned.  Our class includes for example, the algorithms discussed above as
well as multiplier-adjusting types such as the Goldfeld-Quandt-Trotter
approach \cite{GoldQuanTrot66}.  The methods of interest take a
potentially-perturbed Newton step at each iteration so long as the
perturbation is ``not too large'' and the subproblem is solved ``sufficiently
accurately''. The size of the perturbation allowed is simultaneously related
to the parameter $\alpha$ defining the class of methods and the rate of the
asymptotic convergence of the method.  For each method in each
$\alpha$-parametrized class and each $\epsilon \in (0,1)$, we construct a
function with globally $\alpha-$H\"older-continuous Hessian and Lipschitz
continuous gradient for which the method takes precisely
$\lceil\epsilon^{-(2+\alpha)/(1+\alpha)}\rceil$ function evaluations to drive
the gradient norm below $\epsilon$.  As such counts are the same order as the
worst-case upper complexity bound of regularization methods, it follows that
the latter methods are optimal within their respective $\alpha$-class of
methods. As $\alpha$ approaches zero, the worst-case complexity of these
methods approaches that of steepest descent, while for $\alpha=1$, we recover
that of cubic regularization.  We also improve the examples proposed in
\cite{CartGoulToin10a,CartGoulToin11c} in two ways.  The first is that we now
employ objective functions with bounded range, which allows refining the
associated definition of sharp worst-case evaluation complexity bounds, the
second being that the new examples now have finite isolated global minimizers.

The structure of the paper is as follows. Section 2 describes the
parameter-dependent class of methods and objectives of interest; Section 2.1
gives properties of the methods such as their connection to fast asymptotic
rates of convergence while Section 2.2 reviews some well-known examples of
methods covered by our general definition of the class. Section 3 then
introduces two examples of inefficiency of these methods and Section~4
discusses the consequences of these examples regarding the sharpness and
possible optimality of the associated worst-case evaluation complexity
bounds. Further consequences of our results  on the new class proposed by
\cite{CurtRobiSama17c} and \cite{RoyeWrig17} are developed in Section~5 and 6,
respectively. Section 7 draws our conclusions.

{\bf Notation.}  Throughout the paper, $\|\cdot\|$ denotes the Euclidean norm
on $\Re^n$, $I$ the $n\times n$ identity matrix, and $\lambda_{\min}(H)$ and
$\lambda_{\max}(H)$ the left- and right-most eigenvalue of any given symmetric
matrix $H$, respectively.  The condition number of a symmetric positive
definite matrix $M$ is denoted by $\kappa(M) \eqdef
\lambda_{\max}(M)/\lambda_{\min}(M)$. If $M$ is only positive-semidefinite
which we denote by $M \succeq 0$,  and $\lambda_{\min}(M)=0$, then
$\kappa(0)\eqdef +\infty$ unless $M=0$, in which case we set
$\kappa(0)\eqdef 1$. Positive definiteness of $M$ is written as $M \succ 0$.

\numsection{A general parametrized class of methods and objectives}
\label{intro}

Our aim is to minimize a given $C^2$ objective function $f(x)$, $x\in \Re^n$.
We consider methods that generate sequences of iterates $\{x_k\}$ 
for which $\{f(x_k)\}$ is monotonically decreasing, we let
$$
f_k\eqdef f(x_k), \quad g_k\eqdef g(x_k) \tim{and}
H_k\eqdef H(x_k).
$$
where $g(x) = \nabla_x f(x)$ and $H(x) =  \nabla_{xx} f(x)$.

Let $\alpha \in [0,1]$ be a fixed parameter and consider iterative methods 
whose iterations are defined as follows. Given some $x_0\in \Re^n$, let
\beqn{its}
x_{k+1}=x_k+s_k,\quad k\geq 0,
\eeqn 
where $s_k$ satisfies
\beqn{skdef}
(H_k+M_k)s_k= -g_k + r_k
\tim{ with } \|r_k\| \leq \min\left[ \kappa_{rg}\|g_k\|, \kappa_{rs} \|M_ks_k\| \right]
\eeqn
for some residual $r_k$ and constants $\kappa_{rg} \in [0,1)$ and $\kappa_{rs}>0$,
and for some symmetric matrix $M_k$ such that
\beqn{LAMBDAprop}
M_k \succeq 0,
\ms
H_k+ M_k \succeq 0
\eeqn
and
\beqn{lambdacondAL}
\lambda_{\min}(H_k)+\lambda_{\min}(M_k)\leq 
\kappa_{\lambda}\max\left\{|\lambda_{\min}(H_k)|,\|g_k\|^{\frac{\alpha}{1+\alpha}}\right\}
\eeqn
for some $\kappa_{\lambda}>1$  independent of $k$. Without loss of generality,
we assume that $s_k\neq 0$. Furthermore, we require
that no infinite steps are taken, namely 
\beqn{sklbdUPPER}
\|s_k\|\leq \kappa_s
\eeqn
for some $\kappa_s>0$ independent of $k$. \emph{The $\calM.\alpha$ class of
second-order methods consists of all methods whose iterations
satisfy \req{its}--\req{sklbdUPPER}.} The particular choices
$M_k = \lambda_k I$ and $M_k= \lambda_k N_k$ (with $N_k$ symmetric, positive
definite and with bounded condition number) will be of particular interest in
what follows\footnote{Note that \req{lambdacondAL} is slightly more
general than a maybe more natural condition involving $\lambda_{\min} (H_k+M_k)$
instead of $\lambda_{\min} (H_k)+\lambda_{\min}(M_k)$.}. Note that the
definition of $\calM.\alpha$ just introduced generalizes that of M.$\alpha$ in
\cite{CartGoulToin11c}. 

Typically, the expression \req{skdef} for $s_k$ is derived by minimizing 
(possibly approximately) the second-order model
\beqn{model}
m_k(s)=f_k+ g_k^Ts+\half s^T(H_k+\beta_k M_k)s,
\tim{ with }
\beta_k\eqdef\beta_k(s)\geq 0 \tim{ and } \beta_k\leq 1
\eeqn
of $f(x_k+s)$---possibly with an explicit regularizing constraint---with the
aim of obtaining a sufficient decrease of $f$ at the new iterate
$x_{k+1}=x_k+s_k$ compared to $f(x_k)$.  In the definition of an $\calM.\alpha$
method however, the issue of (sufficient) objective-function decrease is not
explicitly addressed/required.  There is no loss of generality in doing so
here since although local refinement of the model may be required to ensure
function decrease, the number of function evaluations to do so (at least for
known methods) does not increase the overall worst-case evaluation complexity
by more than a constant multiple and thus does not affect quantitatively the
worst-case bounds derived; see for example, \cite{CartGoulToin10a,
CartGoulToin11d, GratSartToin08} and also Section 2.2.  Furthermore, the
examples of inefficiency proposed in Section~3 are constructed in such a way
that each iteration of the method automatically provides sufficient decrease
of $f$.

Having defined the classes of methods we shall be concerned with, we
now specify the problem classes that we shall apply the methods in each class
to, in order to demonstrate slow convergence. Given a method in
$\calM.\alpha$, we are interested in minimizing functions $f$ that satisfy
\begin{description}
\item{\framebox[1.4cm]{A.$\alpha$}}
$f:\Re^n\rightarrow \Re$ is twice continuously differentiable and
bounded below, with gradient $g$ being globally Lipschitz continuous on
$\Re^n$ with constant $L_g$, namely, 
\beqn{LipsgNEW}
\|g(x)-g(y)\|\leq L_g\|x-y\|,\tim{for all $x,\,y\in \Re^n$;}
\eeqn
and the Hessian $H$ being globally $\alpha-$H\"older continuous on $\Re^n$ 
with constant $L_{H,\alpha}$, i.e.,
\beqn{LipsHNEW}
 \|H(x)-H(y)\|\leq L_{H,\alpha}\|x-y\|^{\alpha},\tim{for all $x,y \in \Re^n$.}
\eeqn
\hfill$\Box$
\end{description}

\noindent
The case when $\alpha=1$ in A.$\alpha$ corresponds to the Hessian of
$f$ being globally Lipschitz continuous.  Moreover, \req{LipsgNEW} implies
\req{LipsHNEW} when $\alpha=0$, so that the A.$0$ class is that of twice
continuously differentiable functions with globally Lipschitz continuous
gradient. Note also that \req{LipsgNEW} and the existence of $H(x)$ imply that
\beqn{Hbounded}
\|H(x)\| \leq L_g
\eeqn
for all $x \in \Re^n$ \cite[Lemma 1.2.2]{Nest04}, and that every function $f$
satisfying A.$\alpha$ with $\alpha>1$ must be quadratic.  As we will see
below, it turns out that we could weaken the conditions defining A.$\alpha$ by
only requiring \req{LipsgNEW} and \req{LipsHNEW} to hold in an open set
containing all the segments $[x_k,x_k+s_k]$ (the ``path of iterates''), but
these segments of course depend themselves on $f$ and the method applied. 

The next subsection provides some background and justification for the
technical condition \req{lambdacondAL} by relating it to fast rates of
asymptotic convergence, which is a defining feature of second-order
algorithms. In Section~2.2, we then review some methods belonging to
$\calM.\alpha$.

\subsection{Properties of the methods in $\calM.\alpha$}

We first state inclusions properties for $\calM.\alpha$ and A.$\alpha$. 

\llem{MA-inclusion-l}{\mbox{}\\*[-3ex]
\begin{enumerate}
\item Consider a method of $\calM.\alpha_1$  for $\alpha_1 \in [0,1]$ and assume
  that it generates bounded gradients. Then it belongs to $\calM.\alpha_2$ for
  $\alpha_2 \in [0, \alpha_1]$.
\item A.$\alpha_1$ implies A.$\alpha_2$ for $\alpha_2 \in [0, \alpha_1]$, with
  $L_{H,\alpha_2} = \max[ L_{H,\alpha_1}, 2L_g ]$.
\end{enumerate}
}

\proof{
By assumption, $\|g_k\| \leq \kappa_g$ for some $\kappa_g \geq 1$. 
Hence, if $\|g_k\| \geq 1$,
\beqn{gkrel}
\|g_k\|^{\frac{\alpha_1}{1+\alpha_1}}
\leq \kappa_g^{\frac{\alpha_1}{1+\alpha_1}}
\leq \kappa_g
\leq \kappa_g \|g_k\|^{\frac{\alpha_2}{1+\alpha_2}}
\eeqn
for any $\alpha_2 \in [0, \alpha_1]$.  Moreover, \req{gkrel} also holds if
$\|g_k\| \leq 1$, proving the first statement of the lemma.
Now we obtain from \req{Hbounded}, that, if $\|x-y\| > 1$, then
\[
\|H(x) - H(y)\| \leq \|H(x)\|+\|H(y)\| \leq 2L_g \leq 2L_g \|x-y\|^\alpha
\]
for any $\alpha \in [0,1]$. When $\|x-y\|\leq 1$, we may deduce from
\req{LipsHNEW} that, if $\alpha_1 \geq \alpha_2$, then \req{LipsHNEW} with
$\alpha = \alpha_1$ implies \req{LipsHNEW} with $\alpha = \alpha_2$. This
proves the second statement.
} %% epr
\noindent
Observe  if a method is known to be   globally convergent in the sense that
$\|g_k\| \rightarrow 0$ when   $k \rightarrow \infty$, then it obviously
generates bounded gradients and thus the globally convergent methods of
$\calM.\alpha_1$ are included in $\calM.\alpha_2$ ($\alpha_2 \in [0, \alpha_1]$).

We next give a sufficient, more  concise, condition on the
algorithm-generated matrices $M_k$ that implies the bound
\req{lambdacondAL}.

\llem{lemma:Mprop}{Let \req{skdef} and \req{LAMBDAprop} hold.  Assume
  also that the algorithm-generated matrices $M_k$ satisfies
\beqn{Msk}
\lambda_{\min}(M_k) \leq \overline{\kappa}_{\lambda}\|s_k\|^{\alpha},
  \tim{for some $\overline{\kappa}_{\lambda}>1$  and $\alpha \in [0,1]$
       independent of $k$.}
\eeqn
Then \req{lambdacondAL} holds with $\kappa_{\lambda}\eqdef
2\overline{\kappa}_{\lambda}^{\frac{1}{1+\alpha}}(1+\kappa_{rg})$.}

\proof{Clearly, \req{lambdacondAL} holds when $\lambda_{\min}(H_k+M_k)=0$.
When $\lambda_{\min}(H_k+M_k)>0$ and hence $H_k+M_k\succ 0$, \req{skdef}
implies that 
\beqn{skblambda}
\|s_k\|
\leq \frac{\|g_k\|+\|r_k\|}{\lambda_{\min}(H_k+M_k)}
\leq \frac{(1+\kappa_{rg})\|g_k\|}{\lambda_{\min}(H_k)+\lambda_{\min}(M_k)}.
\eeqn
This and \req{Msk} give the inequality
\beqn{Linequation}
\psi(\lambda_{\min}(M_k)) \leq 0
\tim{ with } \psi(\lambda) \eqdef
\lambda^{\frac{1}{\alpha}}(\lambda+ \lambda_{\min}(H_k))
-\overline{\kappa}_{\lambda}^{\frac{1}{\alpha}}(1+\kappa_{rg})\|g_k\|.
\eeqn
Now note that $\psi(0) = \psi(-\lambda_{\min}(H_k)) =
-\overline{\kappa}_{\lambda}^{\frac{1}{\alpha}}(1+\kappa_{rg})\|g_k\|$
and thus
\beqn{phil1}
\psi( \lambda_{1,k} ) <0
\tim{ with }
\lambda_{1,k} = \max\{0,-\lambda_{\min}(H_k)\}.
\eeqn
Moreover, the form of $\psi(\lambda)$ implies that
$\psi(\lambda)$  is strictly increasing
for $\lambda \geq \lambda_{1,k}$.
Define now
\beqn{LV}
\lambda_{2,k}
\eqdef -\lambda_{\min}(H_k)+2\max\left\{|\lambda_{\min}(H_k)|,
\overline{\kappa}_{\lambda}^{\frac{1}{1+\alpha}}(1+\kappa_{rg})^{\frac{\alpha}{1+\alpha}}
\|g_k\|^{\frac{\alpha}{1+\alpha}}\right\}
> \lambda_{1,k}.
\eeqn
Suppose first that $\lambda_{\min}(H_k) < 0$ and
$|\lambda_{\min}(H_k)| \geq \overline{\kappa}_{\lambda}^{\frac{1}{1+\alpha}}
(1+\kappa_{rg})^{\frac{\alpha}{1+\alpha}}\|g_k\|^{\frac{\alpha}{1+\alpha}}$.  Then one verifies that
$\lambda_{2,k} = 3|\lambda_{\min}(H_k)|$ and
\[
\begin{array}{lcl}
\psi(\lambda_{2,k})
& = & (3|\lambda_{\min}(H_k)|)^{\frac{1+\alpha}{\alpha}} -
      (3|\lambda_{\min}(H_k)|)^{\frac{1}{\alpha}}|\lambda_{\min}(H_k)|
-\overline{\kappa}_{\lambda}^{\frac{1}{1+\alpha}}(1+\kappa_{rg})^{\frac{\alpha}{1+\alpha}}
\|g_k\|\\*[1ex]
& = & 2 \cdot
3^{\frac{1}{\alpha}}|\lambda_{\min}(H_k)|^{\frac{1+\alpha}{\alpha}}
-\overline{\kappa}_{\lambda}^{\frac{1}{1+\alpha}}(1+\kappa_{rg})^{\frac{\alpha}{1+\alpha}}\|g_k\|
> 0%\\*[1ex]
%& > & 0.
\end{array}
\]
Suppose now that $\lambda_{\min}(H_k) \geq 0$ and
$|\lambda_{\min}(H_k)| \geq \overline{\kappa}_{\lambda}^{\frac{1}{1+\alpha}}
(1+\kappa_{rg})^{\frac{\alpha}{1+\alpha}}\|g_k\|^{\frac{\alpha}{1+\alpha}}$. Then $\lambda_{2,k}=
\lambda_{\min}(H_k)$ and
\[
\psi(\lambda_{2,k})
= (\lambda_{\min}(H_k))^{\frac{1+\alpha}{\alpha}} +
  (\lambda_{\min}(H_k))^{\frac{1}{\alpha}}|\lambda_{\min}(H_k)|
  -\overline{\kappa}_{\lambda}^{\frac{1}{1+\alpha}}(1+\kappa_{rg})^{\frac{\alpha}{1+\alpha}}\|g_k\|
> 0.
\]
Thus we deduce that $\psi(\lambda_{2,k}) > 0$ whenever $|\lambda_{\min}(H_k)|
\geq \overline{\kappa}_{\lambda}^{\frac{1}{1+\alpha}}
(1+\kappa_{rg})^{\frac{\alpha}{1+\alpha}}\|g_k\|^{\frac{\alpha}{1+\alpha}}$.
Moreover the same inequality obviously holds if
$|\lambda_{\min}(H_k)| < \overline{\kappa}_{\lambda}^{\frac{1}{1+\alpha}}
(1+\kappa_{rg})^{\frac{\alpha}{1+\alpha}}\|g_k\|^{\frac{\alpha}{1+\alpha}}$
because $\psi(\lambda)$ is increasing with $\lambda$. As a consequence,
$\psi(\lambda_{2,k}) > 0$ in all cases. We now combine this inequality,
\req{phil1} and the monotonicity of $\psi(\lambda)$ for $\lambda \geq
\lambda_{1,k}$ to obtain that either 
$\lambda_{\min}(M_k) \leq  \lambda_{1,k} < \lambda_{2,k}$
or $\lambda_{\min}(M_k) \in [\lambda_{1,k},\lambda_{2,k})$ because of of
\req{Linequation}.
Thus $\lambda_{\min}(M_k)\leq \lambda_{2,k}$, which, due to \req{LV} and
$\overline{\kappa}_{\lambda}>1$, implies \req{lambdacondAL}. 
} %epr

\noindent
Thus a method satisfying \req{its}--\req{sklbdUPPER} and \req{Msk} belongs to
$\calM.\alpha$, but not every method in $\calM.\alpha$ needs to satisfy \req{Msk}.
This latter requirement implies the following property regarding
the length of the step generated by methods in $\calM.\alpha$ satisfying \req{Msk}
when applied to functions satisfying A.$\alpha$.

\llem{steppropertyLemma}{Assume that an objective function $f$ satisfying
A.$\alpha$ is minimized by a method satisfying \req{its}, \req{skdef},  
\req{Msk} and such that the conditioning of $M_k$ is bounded in that
$\kappa(M_k)\leq \kappa_\kappa$ for some $\kappa_\kappa \geq 1$. Then
there exists  $\overline{\kappa}_{s,\alpha}>0$
independent of $k$ such that, for $k \geq 0$,
\beqn{sklbd}
\|s_k\|\geq \overline{\kappa}_{s,\alpha} \|g_{k+1}\|^{\frac{1}{1+\alpha}}.
\eeqn
}
\proof{The triangle inequality provides
\beqn{tempg}
\|g_{k+1}\| \leq \| g_{k+1}-(g_k+H_ks_k)\| + \|g_k+H_ks_k\|.
\eeqn
From \req{its}, $g_{k+1}=g(x_k+s_k)$ and Taylor expansion provides 
$g_{k+1}=g_k+ \int_0^1 H(x_k+ \tau s_k )s_k d\tau$.
This and
\req{LipsHNEW} now imply 
\[
\|g_{k+1}-(g_k+H_ks_k)\|\leq \left\|\int_0^1 [ H(x_k+ \tau s_k ) -H(x_k)] d\tau
\right\|\cdot\|s_k\|\leq L_{H,\alpha} (1+\alpha)^{-1} \|s_k\|^{1+\alpha}, 
\]
so that \req{tempg} and  \req{skdef} together give that
\[
\begin{array}{lcl}
\|g_{k+1}\|
%&\leq &L_{H,\alpha}(1+\alpha)^{-1}\|s_k\|^{1+\alpha}+\|M_ks_k\| +\|r_k\| \\*[1ex]
&\leq &L_{H,\alpha}(1+\alpha)^{-1}\|s_k\|^{1+\alpha}+(1+\kappa_{rs})\|M_k\|\,\|s_k\|.
\end{array}
\]
If $M_k \neq 0$, this inequality and the fact that $\kappa(M_k)$ is bounded
then imply that
\[
\|g_{k+1}\|
\leq L_{H,\alpha}(1+\alpha)^{-1}\|s_k\|^{1+\alpha}
     +\kappa(M_k)(1+\kappa_{rs})\lambda_{\min}(M_k)\,\|s_k\|,
\]
while we may ignore the last term on the right-hand side if $M_k=0$.  Hence,
in all cases,
\[
\|g_{k+1}\|
\leq L_{H,\alpha}(1+\alpha)^{-1}\|s_k\|^{1+\alpha}
     +\kappa_\kappa(1+\kappa_{rs})\lambda_{\min}(M_k)\,\|s_k\|,
\]
where we used that $\kappa(M_k)\leq \kappa_\kappa$ by assumption.
This bound and \req{Msk} then imply \req{sklbd} with  
$\overline{\kappa}_{s,\alpha}
\eqdef [L_{H,\alpha}(1+\alpha)^{-1}
         +\kappa_\kappa(1+\kappa_{rs})\overline{\kappa}_{\lambda}]^{-\frac{1}{1+\alpha}}$.
~}  % epr
\noindent
Property \req{sklbd} will be central for proving (in Appendix A2) desirable properties of
a class of methods belonging to $\calM.\alpha$. In addition, we now show that
\req{sklbd} is a necessary condition for fast local convergence of methods of
type \req{skdef}, under reasonable assumptions; fast local rate of convergence
in a neighbourhood of well-behaved minimizers is a ``trademark'' of what is
commonly regarded as second-order methods. 

\llem{lemma:quadratic}{Let $f$ satisfy assumptions A.$\alpha$.  
Apply an algorithm to minimizing $f$ that satisfies
\req{its} and \req{skdef}  and for which
\beqn{multiplierbd}
\|M_k\|\leq \overline{\kappa}_{\lambda}, \tim{
  $k\geq 0$, \quad for some $\overline{\kappa}_{\lambda}>0$ independent
  of $k$.}
\eeqn
Assume also that  convergence at linear or faster than linear rate occurs,
namely,
\beqn{Qconv}
\|g_{k+1}\|\leq \kappa_c \|g_k\|^{1+\alpha}, \quad k\geq 0,
\eeqn
for some $\kappa_c>0$ independent of $k$, with $\kappa_c\in (0,1)$
when $\alpha=0$.  Then \req{sklbd} holds.}

\proof{Let
\beqn{hypo}
0\leq \alpha_k\eqdef \frac{\|s_k\|}{\|g_{k+1}\|^{\frac{1}{1+\alpha}}},\quad k\geq 0.
\eeqn
From \req{Qconv} and the definition of $\alpha_k$ in
\req{hypo}, we have that, for $k \geq 0$,
\[
\begin{array}{lcccl}
  (1-\kappa_{rg})\frac{\|s_k\|}{\alpha_k}
  & \leq & \kappa_{c,\alpha}(1-\kappa_{rg})\|g_k\|
  & \leq & \kappa_{c,\alpha}\|g_k+r_k\| \\*[1ex]
  &   =  & \kappa_{c,\alpha}\|(H_k+M_k)s_k\| 
  & \leq & \kappa_{c,\alpha}\|H_k+M_k\|\cdot \|s_k\|,
\end{array}
\]
where $\kappa_{c,\alpha}\eqdef \kappa_c^{\frac{1}{1+\alpha}}$ and 
where we used \req{skdef} to obtain the first equality. It follows that
\beqn{hl}
\|H_k+M_k\|\geq \frac{(1-\kappa_{rg})}{\alpha_k\kappa_{c,\alpha}},\quad k\geq
0.
\eeqn
The bounds \req{Hbounded} and  \req{multiplierbd} imply that
$\{H_k+M_k \}$ is uniformly bounded above for all $k$,
namely,
\beqn{hl2}
\|H_k+M_k\|\leq \kappa_{hl},\quad k\geq 0,
\eeqn
where $\kappa_{hl}\eqdef L_g+\overline{\kappa}_{\lambda}$.
Now \req{hl} and \req{hl2} give that
$\alpha_k\geq 1/(\kappa_{hl}\kappa_{c,\alpha})>0$, for all $k\geq 0$,
and so it follows from \req{hypo}, that \req{sklbd} holds
with $\overline{\kappa}_{s,\alpha}\eqdef (1-\kappa_{rg})/(\kappa_{c_1}\kappa_{c,\alpha})$.}

\noindent
It is clear from the proof of Lemma \ref{lemma:quadratic} that \req{Qconv} is
only needed asymptotically, that is for all $k$ sufficiently large; for
simplicity, we have assumed it holds globally.

Note that letting $\alpha=1$ in Lemma \ref{lemma:quadratic} provides a
necessary condition for quadratically convergent methods satisfying \req{its},
\req{skdef} and \req{multiplierbd}. Also, similarly to the above proof, one
can show that if superlinear convergence of $\{g_k\}$ to zero occurs, then
\req{sklbd} holds with $\alpha=0$ for all $\overline{\kappa}_{s,\alpha}>0$, or
equivalently, $\|g_{k+1}\|/\|s_k\|\rightarrow 0$, as $k\rightarrow \infty$.

Summarizing, we have shown that \req{sklbd} holds for a
method in $\calM.\alpha$ if \req{Msk} holds and $\kappa(M_k)$ is bounded, or if
linear of faster asymptotic convergence takes place for unit steps.

\subsection{Some examples of methods that belong to the class $\calM.\alpha$}
\label{NtrARC}

Let us now illustrate some of the methods that either by construction or under
certain conditions belong to $\calM.\alpha$.  This list of methods does not
attempt to be exhaustive and other practical methods may be found to belong to
$\calM.\alpha$.

\vspace*{2mm}

\noindent
{\bf Newton's method} \cite{DennSchn83}.\quad 
Newton's method for convex optimization is characterised by finding a
correction $s_k$ that satisfies
$
H_k s_k = - g_k
$
for nonzero $g_k \in {\rm Range}(H_k)$. Letting 
\beqn{Newton}
M_k=0, \ms r_k=0 \tim{and} \beta_k=0
\eeqn
in \req{skdef} and \req{model}, respectively, yields Newton's method. Provided
additionally that both $g_k \in {\rm Range}(H_k)$ and $H_k$ is positive
semi-definite, $s_k$ is a descent direction and \req{LAMBDAprop} holds.  Since
\req{lambdacondAL} is trivially satisfied in this case, it follows that
Newton's method belongs to the class $\calM.\alpha$, for any $\alpha \in [0,1]$,
provided it does not generate infinite steps to violate \req{sklbdUPPER}.  As
Newton's method is commonly embedded within trust-region or regularization
frameworks when applied to nonconvex functions, \req{sklbdUPPER} will in fact,
hold as it is generally enforced for the latter methods.
Note that allowing $\|r_k\| > 0$ subject to the second part of \req{skdef}
then covers inexact variants of Newton's method.

\vspace*{2mm}

\noindent
{\bf Regularization algorithms} \cite{Grie81,
Nest04,CartGoulToin11d}.\quad In these methods, the step $s_k$ from 
the current iterate $x_k$ is computed by (possibly approximately) globally
minimizing the model 
\beqn{r}
m_k(s) = f_k+g_k^Ts+\half s^TH_ks + \frac{\sigma_k}{2+\alpha}\|s\|^{2+\alpha},
\eeqn
where the regularization weight $\sigma_k$ is adjusted to ensure sufficient
decrease of $f$ at $x_k + s_k$.
We assume here that the minimization of \req{r} is carried accurately enough
to ensure that $\nabla_{ss}^² m_k-(s) = H_k+\sigma_k\|s\| I$ is positive
semidefinite, which is always possible because of \cite[Theorem~3.1]{CartGoulToin11}.
The scalar $\alpha$ is the same fixed parameter as in the definition of
A.$\alpha$ and $\calM.\alpha$, so that for each $\alpha \in [0,1]$, we have a
different regularization term and hence what we shall call an
$(2+\alpha)$-regularization method.  For $\alpha=1$, we recover the cubic
regularization (ARC) approach \cite{Grie81,WeisDeufErdm07,NestPoly06,CartGoulToin11, 
CartGoulToin11d}.  For $\alpha=0$, we obtain a quadratic regularization scheme,
reminiscent of the Levenberg-Morrison-Marquardt method \cite{NoceWrig99}.
For these $(2+\alpha)$-regularization methods, we have
\beqn{lbCAL}
\alpha\in [0,1],\quad 
M_k=\sigma_k\|s_k\|^{\alpha} I,
\tim{ and }
\beta_k=\frac{2}{2+\alpha}
\eeqn
in \req{skdef} and \req{model}. If scaling the regularization term is
considered, then the second of these relation is replaced by
$M_k=\sigma_k\|s_k\|^{\alpha} N_k$ for some fixed scaling symmetric positive
definite matrix having a bounded condition number. Note that, by construction,
$\kappa(M_k)=1$. Since $\alpha \geq 0$, we have $0\leq
\beta_k\leq 1$ which is required in \req{model}.  A mechanism of
successful and unsuccessful iterations and $\sigma_k$ adjustments can be
devised similarly to ARC \cite[Alg. 2.1]{CartGoulToin11} in order to deal with steps
$s_k$ that do not give sufficient decrease in the objective. An upper bound on
the number of unsuccessful iterations which is constant multiple of successful
ones can be given under mild assumptions on $f$
\cite[Theorem~2.1]{CartGoulToin11d}. Note that each (successful or unsuccessful)
iteration requires one function- and at most one gradient evaluation.

We now show that for each $\alpha\in [0,1]$, the $(2+\alpha)-$regularization
method based on the model \req{r} satisfies \req{sklbdUPPER} and
\req{lambdacondAL} when applied to $f$ in A.$\alpha$, and so it belongs to
$\calM.\alpha$.

\llem{lemma:regularization}
{Let $f$ satisfy A.$\alpha$ with $\alpha \in   (0,1]$.
Consider minimizing $f$ by applying an $(2+\alpha)$-regularization method
based on the model \req{r}, where the step $s_k$ is chosen as the global
minimizer of the local $\alpha-$model, namely of $m_k(s)$ in \req{model} with
the choice \req{lbCAL}, and where the regularization parameter $\sigma_k$ is
chosen to ensure that
\beqn{sigmamin}
\sigma_k\geq \sigma_{\min},\quad k\geq 0,
\eeqn
for some $\sigma_{\min}>0$ independent of $k$. Then \req{sklbdUPPER} and
\req{Msk} hold,  and so the $(2+\alpha)$-regularization method belongs to
$\calM.\alpha$.
}

\proof{(see Appendix A2 for details)
The same argument that is used in \cite[Lem.2.2]{CartGoulToin11} 
for the $\alpha = 1$ case (see also Appendix A2) provides 
\beqn{reg-sk-upper}
\|s_k\| \leq \max \left \{
\left( \frac{3(2+\alpha) L_g}{4\sigma_k} \right )^{\frac{1}{\alpha}}, 
\left(\frac{3(2+\alpha) \|g_k\|}{\sigma_k}\right)^{\frac{1}{1+\alpha}} 
\right \}, \quad k\geq 0,
\eeqn
so long as A.$\alpha$ holds, which together with \req{sigmamin}, implies
\beqn{sigmaCAUCHY}
\|s_k\| \leq \max \left \{
\left( \frac{3(2+\alpha) L_g}{4\sigma_{\min}} \right )^{\frac{1}{\alpha}}, 
\left(\frac{3(2+\alpha) \|g_k\|}{\sigma_{\min}}\right)^{\frac{1}{1+\alpha}} 
\right \}, \quad k\geq 0.
\eeqn
The assumptions A.$\alpha$, that
%we employ the true Hessians rather than approximations and that
the model is minimized globally imply that the 
$\alpha \leq 1$ analog of 
\cite[Corollary 2.6]{CartGoulToin11} holds, which gives
$\|g_k\|\rightarrow 0$ as $k\rightarrow \infty$, and so $\{\|g_k\|\}$,
$k\geq 0$, is bounded above. The bound
\req{sklbdUPPER} now follows from \req{sigmaCAUCHY}.

Using the same techniques as in \cite[Lemma 5.2]{CartGoulToin11} that applies when
$f$ satisfies A.$1$, it is easy to show for the more general A.$\alpha$ case
that $\sigma_k \leq c_{\sigma}\max(\sigma_0,L_{H,\alpha})$ for all $k$,  where $c_{\sigma}$
is a constant depending solely on $\alpha$ and algorithm parameters. It then
follows from \req{lbCAL} that \req{Msk} holds and therefore that the 
$(2+\alpha)$-regularization method belongs to $\calM.\alpha$ for $\alpha \in (0,1]$. 
}

\noindent
We cannot extend this result to the $\alpha = 0$ case unless we also assume
that $H_k$ is positive semi-definite. If this is the case, further examination
of the proof of \cite[Lem.2.2]{CartGoulToin11} allows us to remove the
first term in the max in \req{sigmaCAUCHY}, and the remainder of the proof is
valid.

We note that bounding the regularization parameter $\sigma_k$ away from zero
in \req{sigmamin} appears crucial when establishing the bounds
\req{sklbdUPPER} and \req{lambdacondAL}.  Requiring \req{sigmamin} implies
that the Newton step is always perturbed, but does not prevent local quadratic
convergence of ARC \cite{CartGoulToin11d}.

\vspace*{2mm}

\noindent
{\bf Goldfeld-Quandt-Trotter-type (GQT) methods} \cite{GoldQuanTrot66}.\quad
Let $\alpha \in (0,1]$. These algorithms set $M_k = \lambda_k I$, where
\beqn{gqtCONSTR}
\lambda_k=\left\{
\begin{array}{ll}
0,&\tim{when $\lambda_{\min}(H_k)\geq \omega_k\|g_k\|^{\frac{\alpha}{1+\alpha}}$;}\\
-\lambda_{\min}(H_k)+\omega_k\|g_k\|^{\frac{\alpha}{1+\alpha}},&
  \tim{otherwise,}
\end{array}
  \right.
  \eeqn
in \req{skdef}, where $\omega_k>0$ is a parameter that is adjusted
so as to ensure sufficient objective decrease.  (Observe that replacing
$\frac{\alpha}{1+\alpha}$ by $1$ in the exponent of $\|g_k\|$ in
\req{gqtCONSTR} recovers the original method of Goldfeld et
al. \cite{GoldQuanTrot66}.)
It is straightforward to check that \req{LAMBDAprop} holds for the choice
\req{gqtCONSTR}. Thus the GQT approach takes the pure Newton step whenever the
Hessian is locally sufficiently positive definite, and a suitable
regularization of this step otherwise.  The parameter $\omega_k$ is increased
by a factor, say $\gamma_1>1$, and $x_{k+1}$ left as $x_k$ whenever the step
$s_k$ does not give sufficient decrease in $f$ (i.e., iteration $k$ is
unsuccessful), namely when
\beqn{rho}
\rho_k\eqdef \frac{f_k-f(x_k+s_k)}{f_k-m_k(s_k)}\leq \eta_1,
\eeqn
where $\eta_1\in (0,1)$ and 
\beqn{qmod}
m_k(s) = f_k+g_k^Ts+\half s^TH_ks
\eeqn
is the model \req{model} with $\beta_k=0$.  If $\rho_k>\eta_1$, then
$\omega_{k+1}\leq \omega_k$ and $x_{k+1}$ is constructed as in
\req{its}.
Note that the choice \req{gqtCONSTR} implies that \req{lambdacondAL}
holds, provided $\omega_k$ is uniformly bounded above. We show that the
latter, as well as \req{sklbdUPPER}, hold for functions in A.$\alpha$.

\llem{lemma:GQT}{Let $f$ satisfy A.$\alpha$ with $\alpha \in  (0,1]$.
Consider minimizing $f$ by applying a GQT method that sets $\lambda_k$
in \req{skdef} according to \req{gqtCONSTR}, measures progress
according to  \req{rho}, and chooses the parameter $\omega_k$ and the residual
$r_k$ to satisfy, for $k\geq 0$,
\beqn{Rmin}
\omega_k \geq \omega_{\min} \quad k\geq 0.
\tim{ and }
r_k^Ts_k \leq 0.
\eeqn
Then \req{sklbdUPPER} and \req{lambdacondAL} hold,  and so the GQT
method belongs to $\calM.\alpha$.}
  
Note that the second part of \req{Rmin} merely requires that
$s_k$ is not longer that the line minimum of the regularized model along the
direction $s_k$, that is $1 \leq {\rm arg}\min_{\tau \geq 0}m_k(\tau s_k)$.
  
\proof{
Let us first show \req{sklbdUPPER}. Since $\omega_k>0$, and
$g_k+r_k\neq 0$ until termination, the choice of $\lambda_k$ in 
\req{gqtCONSTR} implies that $\lambda_k+\lambda_{\min}(H_k)>0$, for
all $k$, and so \req{skdef} provides
\beqn{GQT-step}
s_k=-(H_k+\lambda_k I)^{-1}(g_k+r_k),
\eeqn
and hence,
\beqn{skdef2}
\|s_k\|
\leq \|(H_k+\lambda_k I)^{-1}\|\cdot\|g_k+r_k\|
=\frac{(1+\kappa_{rg})\|g_k||}{\lambda_k+\lambda_{\min}(H_k)},\quad k\geq 0.
\eeqn
It follows from \req{gqtCONSTR} and the first part of \req{Rmin} that, for all
$k \geq 0$,  
\beqn{GQT-b1}
\lambda_k+\lambda_{\min}(H_k) 
\geq \omega_k\|g_k\|^{\frac{\alpha}{1+\alpha}}
\geq \omega_{\min}\|g_k\|^{\frac{\alpha}{1+\alpha}},
\eeqn
This and \req{skdef2} further give
\beqn{skdef3}
\|s_k\|\leq \frac{(1+\kappa_{rg})\|g_k\|^{\frac{1}{1+\alpha}}}{\omega_{\min}},\quad k\geq 0.
\eeqn
As global convergence assumptions are satisfied when $f$ in
A.$\alpha$ \cite{ConnGoulToin00, GoldQuanTrot66}, we have $\|g_k\|\rightarrow 0$ as
$k\rightarrow \infty$ (in fact, we only need  the gradients
$\{g_k\}$ to be bounded). Thus \req{skdef3} implies~\req{sklbdUPPER}.

Due to \req{gqtCONSTR}, \req{lambdacondAL} holds if we show that
$\{\omega_k\}$ is uniformly bounded above. For this, we first need
to estimate the model decrease. Taking the inner product of
\req{skdef} with $s_k$, we obtain that
$$
-g_k^Ts_k
=s_k^TH_ks_k+\lambda_k \|s_k\|^2 - r_k^Ts_k.
$$
Substituting this into the model decrease, we deduce also from
\req{model} with $\beta_k=0$ that
$$
f_k-m_k(s_k)
=-g_k^Ts_k-\half s_k^TH_ks_
k=\half s_k^TH_ks_k+\lambda_k\|s_k\|^2- r_k^Ts_k
\geq \left(\half\lambda_{\min}(H_k)+\lambda_k\right)\|s_k\|^2.
$$
where we used the second part of \req{Rmin} to obtain the last inequality.
It is straightforward to check that this and \req{GQT-b1} now imply
\beqn{QGTmodeldecrease}
f_k-m_k(s_k)
\geq \half \omega_k\|g_k\|^{\frac{\alpha}{1+\alpha}}\cdot\|s_k\|^2.
\eeqn
We show next that iteration $k$ is successful for $\omega_k$ sufficiently
large.  From \req{rho} and second-order Taylor expansion of $f(x_k+s_k)$, we
deduce
$$
|\rho_k-1|
=\left|\frac{f(x_k+s_k)-m_k(s_k)}{f_k-m_k(s_k)}\right|
\leq \frac{|H_k-H(\xi_k)|\cdot\|s_k\|^2}{2(f_k-m_k(s_k))}
\leq \frac{L_{H,\alpha}\|s_k\|^{2+\alpha}}{2(f_k-m_k(s_k))}.
$$
This and \req{QGTmodeldecrease} now give
\beqn{rho1}
|\rho_k-1|\leq
  \frac{L_{H,\alpha}\|s_k\|^{\alpha}}{\omega_k\|g_k\|^{\frac{\alpha}{1+\alpha}}}\leq
  \frac{L_{H,\alpha}}{\omega_{\min}^{\alpha}\omega_k},
\eeqn
where to obtain the last inequality, we used \req{skdef3}. Due to
\req{rho}, iteration $k$ is successful when $|\rho_k-1|\leq 1-\eta_1$,
which from \req{rho1} is guaranteed to hold whenever $\omega_k\geq
\frac{L_{H,\alpha}}{\omega_{\min}^{\alpha}(1-\eta_1)}$. As on each
successful iteration we set $\omega_{k+1}\leq \omega_k$, it follows that 
\beqn{Rbd}
\omega_k\leq \overline{\omega}\eqdef
\max\left\{
    \omega_0,\frac{\gamma_1L_{H,\alpha}}{\omega_{\min}^{\alpha}(1-\eta_1)}
\right\},\quad
k\geq 0,
\eeqn
where the $\max$ term addresses the situation at the starting point and the
$\gamma_1$ factor is included in case an iteration was unsuccessful and close
to the bound. This concludes proving \req{lambdacondAL}.
~} %epr

\vspace*{2mm}

\noindent
{\bf Trust-region algorithms} \cite{ConnGoulToin00}.\quad 
These methods compute the correction $s_k$ as the global solution of the
subproblem
\beqn{trs}
\tim{minimize} f_k+g_k^Ts+\half s^TH_ks \tim{subject to}  \|s\|\leq \Delta_k,
\eeqn
where $\Delta_k$ is an evolving trust-region radius that is chosen to
ensure sufficient decrease of $f$ at $x_k + s_k$.
The resulting global minimizer 
satisfies \req{skdef}--\req{LAMBDAprop} \cite[Corollary 7.2.2]{ConnGoulToin00}
with $M_k = \lambda_k I$ (or $M_k = \lambda_k N_k$ if scaling is considered) and
$r_k=0$. The scalar $\lambda_k$ is the Lagrange multiplier of the
trust-region constraint, satisfies 
\beqn{trust-region}
\lambda_k\geq \max\{0,-\lambda_{\min}(H_k)\}
\eeqn
and is such that $\lambda_k=0$ whenever $\|s_k\|<\Delta_k$ (and then, $s_k$ is
the Newton step) or calculated using \req{skdef} to ensure that
$\|s_k\|=\Delta_k$. The scalar $\beta_k=0$ in \req{model}.  The iterates are
defined by \req{its} whenever sufficient progress can be made in some relative
function decrease (so-called {\em successful iterations}), and they remain
unchanged otherwise ({\em unsuccessful iterations}) while $\Delta_k$ is
adjusted to improve the model (decreased on unsuccessful iterations, possibly
increased on successful ones).
The total number of unsuccessful iterations is
bounded above by a constant multiple of the successful ones plus a
(negligible) term in $\log \epsilon$ \cite[page
23]{GratSartToin08} provided $\Delta_k$ is not increased too fast on
successful iterations. One successful iteration requires one gradient and one
function evaluation while an unsuccessful one only evaluates the objective.

The property \req{sklbdUPPER} of $\calM.\alpha$ methods can be easily shown for
trust-region methods, see Lemma \ref{lemma:TR}.  It is unclear however,
whether conditions \req{lambdacondAL} or \req{Msk} can be guaranteed in
general for functions in A.$\alpha$.  The next lemma gives conditions ensuring
a uniform upper bound on the multiplier $\lambda_k$, which still falls short
of \req{lambdacondAL} in general.

\llem{lemma:TR}{Let $f$ satisfy assumptions A.$0$.
Consider minimizing $f$ by applying a trust-region 
method as described in \cite[Algorithm 6.1.1]{ConnGoulToin00}, where
the trust-region subproblem is minimized globally to compute
$s_k$ and where the trust-region radius is chosen to ensure that
\beqn{TRradius}
\Delta_k\leq \Delta_{\max},\quad k\geq 0,
\eeqn
for some $\Delta_{\max}>0$. Then \req{sklbdUPPER} holds.
Additionally, if
\beqn{decreasingG}
\|g_{k+1}\|\leq \|g_k\|, \tim{for all $k$
    sufficiently large,}
\eeqn
then $\lambda_k\leq \lambda_{\max}$ for all $k$ and some $\lambda_{\max}>0$,
and $\lambda_{\min}(M_k)$ is bounded.}
\proof{Consider the basic trust-region algorithm as described in
  \cite[Algorithm 6.1.1]{ConnGoulToin00}, using the same notation.
Since the global minimizer $s_k$ of the trust-region subproblem is feasible
with respect to the trust-region constraint, we have $\|s_k\|\leq \Delta_k$,
and so \req{sklbdUPPER} follows trivially from \req{TRradius}.

Clearly, the upper bound on $\lambda_k$ holds whenever $\lambda_k=0$ 
or $\lambda_k=-\lambda_{\min}(H_k)\leq  L_g$.  Thus
it is sufficient to consider the case when $\lambda_k>0$ and
$H_k+\lambda_k I\succ 0$. The first condition implies that the trust-region 
constraint is active, namely $\|s_k\|=\Delta_k$ \cite[Corollary 7.2.2]{ConnGoulToin00}.
The second condition together with \req{skdef}  implies, as in the proof
of Lemma \ref{lemma:Mprop}, that \req{skblambda} holds. Thus we deduce
$$
\Delta_k\leq \frac{\|g_k\|}{\lambda_k+\lambda_{\min}(H_k)},
$$
or equivalently,
\beqn{wanted}
\lambda_k\leq \frac{\|g_k\|}{\Delta_k}-\lambda_{\min}(H_k)\leq 
 \frac{\|g_k\|}{\Delta_k}+L_g,\quad k\geq 0.
\eeqn
It remains to show that 
\beqn{toprove}
\tim{$\{\|g_k\|/\Delta_k\}$   is bounded above independently of $k$.}
\eeqn
By \cite[Theorem 6.4.2]{ConnGoulToin00},  we have that there exists $c\in (0,1)$ such 
that the implication holds 
\beqn{Dsucc}
\Delta_k\leq c\|g_k\| \quad \Longrightarrow\quad \Delta_{k+1}\geq \Delta_k, 
 \tim{i.e., $k$ is successful.} 
\eeqn
(Observe that the Cauchy model decrease condition \cite[Theorem 6.3.3]{ConnGoulToin00}
 is sufficient to obtain the above implication.)  
Let $\gamma_1\in (0,1)$ denote the largest factor we allow $\Delta_k$ to 
be decreased by (during unsuccessful iterations). Using a similar argument to
that of \cite[Theorem  6.4.3]{ConnGoulToin00}, we let $k\geq k_0$ be the first iterate
such that
\beqn{lowDbound}
\Delta_{k+1}< c\gamma_1\|g_{k+1}\|,
\eeqn
where $k_0$ is the iteration from which onwards \req{decreasingG}
holds. Then since $\Delta_{k+1}\geq \gamma_1\Delta_k$
and from \req{decreasingG} we have that $\Delta_k<c\|g_k\|$. 
This and \req{Dsucc} give 
$$\Delta_{k+1}\geq \Delta_k\geq c\gamma_1\|g_k\|\geq  c\gamma_1\|g_{k+1}\|,$$
where to obtain the second and third inequalities, we used the
hypothesis and \req{decreasingG}, respectively.
We have reached a contradiction with our assumption that $k+1$ is the first 
iteration greater than $k_0$
such that \req{lowDbound} holds.  Hence there is no such $k$ and we deduce that
\beqn{tempDelta}
\Delta_k\geq \min\left\{\Delta_{k_0}, c\gamma_1\|g_k\|\right\} 
 \tim{for all $k\geq k_0$.}
\eeqn 
Note that since $g_k$ remains unchanged on unsuccessful
iterations,  \req{decreasingG} trivially holds on such iterations.
Since the assumptions of \cite[Theorem 6.4.6]{ConnGoulToin00} are satisfied,
we have that $\|g_k\|\rightarrow 0$, as $k\rightarrow \infty$. This
and \req{tempDelta} imply \req{toprove}.  The desired conclusion then follows
from \req{wanted}.}

\noindent
Note that if \req{Qconv} holds for some $\alpha\in [0,1]$, then
\req{decreasingG} is satisfied, and so Lemma \ref{lemma:TR} shows that if
\req{Qconv} holds, then \req{multiplierbd} is satisfied. It follows from Lemma
\ref{lemma:quadratic} that fast convergence of trust-region methods for
functions in A.$\alpha$ alone is sufficient to ensure \req{sklbd}, which in
turn is connected to our definition of the class $\calM.\alpha$.  However, the
properties of the multipliers (in the sense of \req{lambdacondAL} for any
$\alpha\in [0,1]$ or even \req{sklbd}) remain unclear in the absence of fast
convergence of the method.
Based on our experience, we are inclined to believe that generally,
the multipliers $\lambda_k$ are at best guaranteed to be uniformly bounded
above, even for specialized, potentially computationally expensive, rules of
choosing the trust-region radius.

As the Newton step is taken in the trust-region framework satisfying
\req{skdef} whenever it is within the trust region and gives sufficient
decrease in the presence of local convexity, the A.$1$- (hence A.$\alpha$-)
example of inefficient behaviour for Newton's method of worst-case evaluation
complexity precisely $\epsilon^{-2}$ can be shown to apply also to
trust-region methods \cite{CartGoulToin10a} (see also \cite{GratSartToin08}).

\vspace*{2mm}

\noindent
{\bf Linesearch methods} \cite{DennSchn83,NoceWrig99}.\quad
We finally consider methods using a linesearch to control
improvement in the objective at each step. Such methods
compute $x_{k+1}=x_k+s_k$, $k \geq 0$, where $s_k$ is defined via \req{skdef}
in which $M_k$ is chosen so that $H_k + M_k$, the Hessian of the selected 
quadratic model $m_k(s)$, is ``sufficiently'' positive definite, and
$r_k = (1-\mu_k)g_k$, yielding a stepsize $\mu_k \in [1-\kappa_{rg},1]$ which
is calculated so as to decrease $f$ (the {\em linesearch}); this is always
possible for sufficiently small $\mu_k$ (and hence sufficiently small
$\kappa_{rg}$.) The precise definition of ''sufficient decrease'' depends on
the particular linesearch scheme being considered, but we assume here that
\[
\mu_k = 1 \tim{ is acceptable whenever } m_k(s_k) = f(x_k+s_k).
\]
In other words, we require the unit step to be acceptable when the model and
the true objective function match at the trial point.  Because the minimization
of the quadratic model along the step always ensure that $m_k(s_k)= f(x_k) +
\half g_ks_k$, the above condition says that $s_k$ must be acceptable with
$\mu_k = 1$ whenever $f(x_k+s_k) = f(x_k) + \half g_ks_k$.
This is for instance the case for the Armijo and Goldstein linesearch
conditions\footnote{With reasonable algorithmic constants, see Appendix A1.},
two standard linesearch techniques.  As a consequence, the corresponding
linesearch variants of Newton's method and of the $(2+\alpha)$-regularization
methods also belong to $\calM.\alpha$ (with $\beta_k=1$ for all $k$), and the list
is not exhaustive. Note that linesearch methods where the search direction is
computed inexactly are also covered by setting $r_k = g_k - \mu_k(g_k + w_k)$
for some ``error vector'' $w_k$, provided the second part of \req{skdef} still
holds.

\numsection{Examples of inefficient behaviour}
\label{OneDSection}

After reviewing the methods in $\calM.\alpha$, we now turn to showing they can
converge slowly when applied to specific functions with fixed
range\footnote{At variancewith the examples proposed in
\cite{CartGoulToin10a,CartGoulToin11c}.}  and the relevant degree of smoothness.

\subsection{General methods in $\calM.\alpha$}\label{slowMa-s} 

Let $\alpha\in [0,1]$ and $\epsilon \in (0,1)$ be given and consider an
arbitrary method in $\calM.\alpha$. Our intent is now to construct a univariate
function $f^{\calM.\alpha}_\epsilon(x)$ satisfying A.$\alpha$ such that
\beqn{ex-fbounds}
f^{\calM.\alpha}_\epsilon(0) = 1,
\ms
f^{\calM.\alpha}_\epsilon(x) \in [ a, b ]  \tim{for} x \geq 0,
\eeqn
for some constants $a \leq b $ independent of $\epsilon$ and $\alpha$, and
such that the method will terminate in exactly
\beqn{ex-nespalphadef}
k_{\epsilon,\alpha} = \left\lceil\epsilon^{-\frac{2+\alpha}{1+\alpha}}\right\rceil
\eeqn
iterations (and evaluations of $f$, $g$ and $H$).

We start by defining the sequences ${f_k}$, ${g_k}$ and ${H_k}$ for
$ k = 0, \ldots, k_{\epsilon,\alpha}$ by
\beqn{ex-fkgkHk}
f_k  = 1 -  \half k \epsilon^{\frac{2+\alpha}{1+\alpha}} ,
\ms
g_k = - 2 \,\epsilon \, f_k
\tim{and}
H_k = 4 \,\epsilon^{\frac{\alpha}{1+\alpha}}\, f_k^2.
\eeqn
They are intended to specify the
objective function, gradient and Hessian values at successive iterates
generated by the chosen method in $\calM.\alpha$, according to \req{its}
and \req{skdef} for some choice of multipliers
$\{\lambda_k\} = \{M_k\} = \{\lambda_{\min}(M_k)\}$
satisfying \req{LAMBDAprop} and \req{lambdacondAL}.  In other words,
we impose that $f_k = f^{\calM.\alpha}_\epsilon(x_k)$,
$g_k = \nabla f^{\calM.\alpha}_\epsilon(x_k)$ and $H_k = \nabla^2 f^{\calM.\alpha}_\epsilon(x_k)$
for $k \in \calK \eqdef \iibe{0}{k_{\epsilon,\alpha}}$.
Note that ${f_k}$, ${|g_k|}$ and ${H_k}$ are monotonically decreasing and
that, using \req{ex-nespalphadef},
\beqn{ex-fkbounds}
f_k \in [\half, 1] \tim{ for } k \in \calK.
\eeqn
In addition, \req{LAMBDAprop} and \req{lambdacondAL} impose that, for $k \in \calK$,
\[
0 \leq \lambda_k + 4\epsilon^{\frac{\alpha}{1+\alpha}} f_k^2
\leq \kappa_\lambda \max[ 4\epsilon^{\frac{\alpha}{1+\alpha}} f_k^2,
  (2\epsilon f_k)^{\frac{\alpha}{1+\alpha}} ]
= 4\kappa_\lambda \epsilon^{\frac{\alpha}{1+\alpha}} f_k^2.
\]
yielding that
\beqn{ex-lambdabound}
\lambda_k \in \left[0, 4(\kappa_\lambda-1) \epsilon^{\frac{\alpha}{1+\alpha}}f_k^2\right],
\eeqn
As a consequence, we obtain, using both parts of \req{skdef},
that, for $k \in \calK$,
\beqn{ex-skdef}
s_k = \theta_k \frac{\epsilon^{\frac{1}{1+\alpha}}}{2f_k}
\tim{ for some } \theta_k \in \left[\frac{1-\kappa_{rg}}{\kappa_\lambda},1+\kappa_{rg}\right].
\eeqn
Note that our construction imposes
\comment{%%%%%%%%%%%%%%%%%%%%%%%%%%%%%%%%%%%%%%%%%%%%%%%%%%%%%%%%
  a sufficient decrease property in that
\[
f^{\calM.\alpha}_\epsilon(x_k+s_k)
= f_{k+1}
= f_k - \frac{1}{2} \epsilon^{\frac{2+\alpha}{1+\alpha}}
= f^{\calM.\alpha}_\epsilon(x_k) + \frac{1}{2\theta_k}g_k^Ts_k
\leq f^{\calM.\alpha}_\epsilon(x_k) + \frac{1}{2(1+\kappa_{rg})} g_k^Ts_k,
\]
which can be interpreted as an Armijo-like descent condition.
}%%%%%%%%%%%%%%%%%%%%%%%%%%%%%%%%%%%%%%%%%%%%%%%%%%%%%%%%%%%%%%%%%
that
\beqn{ex-d1}
\begin{array}{lcl}
m_k(s_k)
&   =  & f_k + g_ks_k + \half g_ks_k +\half s_k(H_k+\beta_k\lambda_k)s_k \\
&   =  & f_k + g_ks_k + \half s_k[-g_k+r_k + (\beta_k-1)\lambda_ks_k] \\
& \geq & f_k - \half |g_k|s_k -\half \kappa_{rg}|g_k|s_k
        +\half \theta_k^2 (\kappa_\lambda-1)(\beta_k-1)\epsilon^{\frac{2+\alpha}{1+\alpha}} \\
& \geq & f_k - \half \theta_k\epsilon^{\frac{2+\alpha}{1+\alpha}}
                  [1+\kappa_{rg}+\theta_k(1-\beta_k)(\kappa_\lambda-1)]\\
& \geq &  f_k - \half \epsilon^{\frac{2+\alpha}{1+\alpha}}(1+\kappa_{rg})^2
                  [1 + (1-\beta_k)(\kappa_\lambda-1)]\\
& \geq &  f_k - \half \epsilon^{\frac{2+\alpha}{1+\alpha}}(1+\kappa_{rg})^2\kappa_\lambda
\end{array}                  
\eeqn
where we have used \req{skdef}, \req{ex-fkgkHk}, \req{ex-skdef},
\req{ex-lambdabound} and $\beta_k \leq 1$.
Hence, again taking  \req{ex-fkgkHk} into account,
\beqn{ex-decrease}
\frac{f_k-f_{k+1}}{f_k-m_k(s_k)}
% = \frac{\half \epsilon^{\frac{2+\alpha}{1+\alpha}}}{\half |g_k|s_k (1-\kappa_{rg})}
\geq \frac{\half \epsilon^{\frac{2+\alpha}{1+\alpha}}}{\half \epsilon^{\frac{2+\alpha}{1+\alpha}}\kappa_\lambda(1+\kappa_{rg})^2}
= \frac{1}{(1+\kappa_{rg})^2\kappa_\lambda}
\in (0,1),
\eeqn
and sufficient decrease of the objective function automatically follows.
Moreover, given \req{ex-fkbounds}, we deduce from \req{ex-skdef} that $|s_k|
\leq 1$ for $k\in \calK$ and \req{sklbdUPPER} holds with $\kappa_s=1$, as
requested  for a method in $\calM.\alpha$.
It also follows from \req{its} and \req{ex-skdef} that, if $x_0=0$,
\beqn{sk1}
s_k>0 \tim{and} x_k=\sum_{i=0}^{k-1}s_i, \quad k = 0, \ldots, k_{\epsilon,\alpha}.
\eeqn
We therefore conclude that the sequences
$\{f_k\}_{k=0}^{k_{\epsilon,\alpha}}$, $\{g_k\}_{k=0}^{k_{\epsilon,\alpha}}$,
$\{H_k\}_{k=0}^{k_{\epsilon,\alpha}}$, $\{\lambda_k\}_{k=0}^{k_{\epsilon,\alpha}-1}$ and
$\{s_k\}_{k=0}^{k_{\epsilon,\alpha}-1}$ can be viewed as produced by our
chosen method in $\calM.\alpha$, and, from \req{ex-fkgkHk}, that termination
occurs precisely for $k = k_{\epsilon,\alpha}$, as desired.

We now construct the function $f^{\calM.\alpha}_\epsilon(x)$ for
$x \in [0,x_{k_{\epsilon,\alpha}}]$ using Hermite interpolation. We set
\beqn{fM1}
f^{\calM.\alpha}_\epsilon(x)=  p_k(x-x_k)+f_{k+1}
\tim{for $x\in[x_k,x_{k+1}]$ and $k = 0,   \ldots, k_{\epsilon,\alpha}-1$,}
\eeqn
where $p_k$ is the polynomial
$$
p_k(s)=c_{0,k}+c_{1,k}s+c_{2,k}s^2+c_{3,k}s^3+c_{4,k}s^4+c_{5,k}s^5,
$$
with coefficients defined by the interpolation conditions
\beqn{interpolation2}
\begin{array}{l}
p_k(0)=f_k-f_{k+1},\quad
p_k(s_k)=0;\\[1.5ex]
p_k^{\prime}(0)=g_k,\quad p_k^{\prime}(s_k)=g_{k+1};\\[1.5ex]
p_k^{''}(0)=H_k,\quad p_k^{''}(s_k)=H_{k+1},
\end{array}
\eeqn
where $s_k$ is defined in \req{ex-skdef}.
These conditions yield the following values for the coefficients
\beqn{ex-c0c1}
c_{0,k}=f_k-f_{k+1},\quad
c_{1,k}=g_k,\quad c_{2,k}=\half H_k;
\eeqn
with the remaining coefficients satisfying
$$
\left(
\begin{array}{lcr}
s_k^3&s_k^4&s_k^5\\
3s_k^2&4s_k^3&5s_k^4\\
6s_k&12s_k^2&20s_k^3
\end{array}
\right)
\left(
\begin{array}{c}
c_{3,k}\\
c_{4,k}\\
c_{5,k}
\end{array}
\right)=
\left(
\begin{array}{c}
\Delta f_k-g_ks_k-\half s_k^TH_ks_k\\
\Delta g_k-H_ks_k\\
\Delta H_k
\end{array}
\right), 
$$
where
$$
\Delta f_k=f_{k+1}-f_k,\quad \Delta g_k=g_{k+1}-g_k \tim{and} \Delta H_k=H_{k+1}-H_k.
$$
Hence we obtain after elementary calculations that
\beqn{ex-expck}
\begin{array}{c}
c_{3,k}
=10\bigfrac{\Delta f_k}{s_k^3}-4\bigfrac{\Delta
  g_k}{s_k^2}+\bigfrac{\Delta
  H_k}{2s_k}-10\bigfrac{g_k}{s_k^2}-\bigfrac{H_k}{s_k};\\[2ex]
%= 10\bigfrac{\Delta f_k}{s_k^3}-4\bigfrac{\Delta
%  g_k}{s_k^2}+\bigfrac{\Delta H_k}{2s_k}-9\bigfrac{g_k}{s_k^2}
%   +\bigfrac{\lambda_k}{s_k};\\[2ex]
%= \bigfrac{1}{s_k^2}\left[\epsilon f_k \left( 20 - \bigfrac{10}{\theta_k}-2\theta_k\right)
%  -\epsilon^{\frac{3+2\alpha}{1+\alpha}}(4 + \pi_k)\right];\\[2ex]
c_{4,k}
=-15\bigfrac{\Delta f_k}{s_k^4}+7\bigfrac{\Delta
  g_k}{s_k^3}-\bigfrac{\Delta H_k}{s_k^2}+15\bigfrac{g_k}{s_k^3}+\bigfrac{H_k}{2s_k^2};\\[2ex]
%=-15\bigfrac{\Delta f_k}{s_k^4}+7\bigfrac{\Delta
%  g_k}{s_k^3}-\bigfrac{\Delta H_k}{s_k^2}+\frac{29}{2}\cdot\bigfrac{g_k}{s_k^3}
%    -\bigfrac{\lambda_k}{2s_k^2};\\[2ex]
%= \bigfrac{1}{s_k^3}\left[\epsilon f_k \left(\bigfrac{15}{\theta_k}-30 + \theta_k\right)
%   + \epsilon^{\frac{3+2\alpha}{1+\alpha}}(7 + 2\pi_k)\right];\\[2ex]
c_{5,k}=6\bigfrac{\Delta f_k}{s_k^5}-3\bigfrac{\Delta
  g_k}{s_k^4}+\bigfrac{\Delta H_k}{2s_k^3}-6\bigfrac{g_k}{s_k^4};\\[2ex]
%= \bigfrac{1}{s_k^4}\left[\epsilon f_k \left(12 - \bigfrac{6}{\theta_k}\right)
%  -\epsilon^{\frac{3+2\alpha}{1+\alpha}}(3 + \pi_k)\right],
\end{array}
\eeqn

\begin{figure}[htbp] % produced by exalpha.m
\begin{center}
\vspace*{1.5mm}
\includegraphics[height=3.5cm]{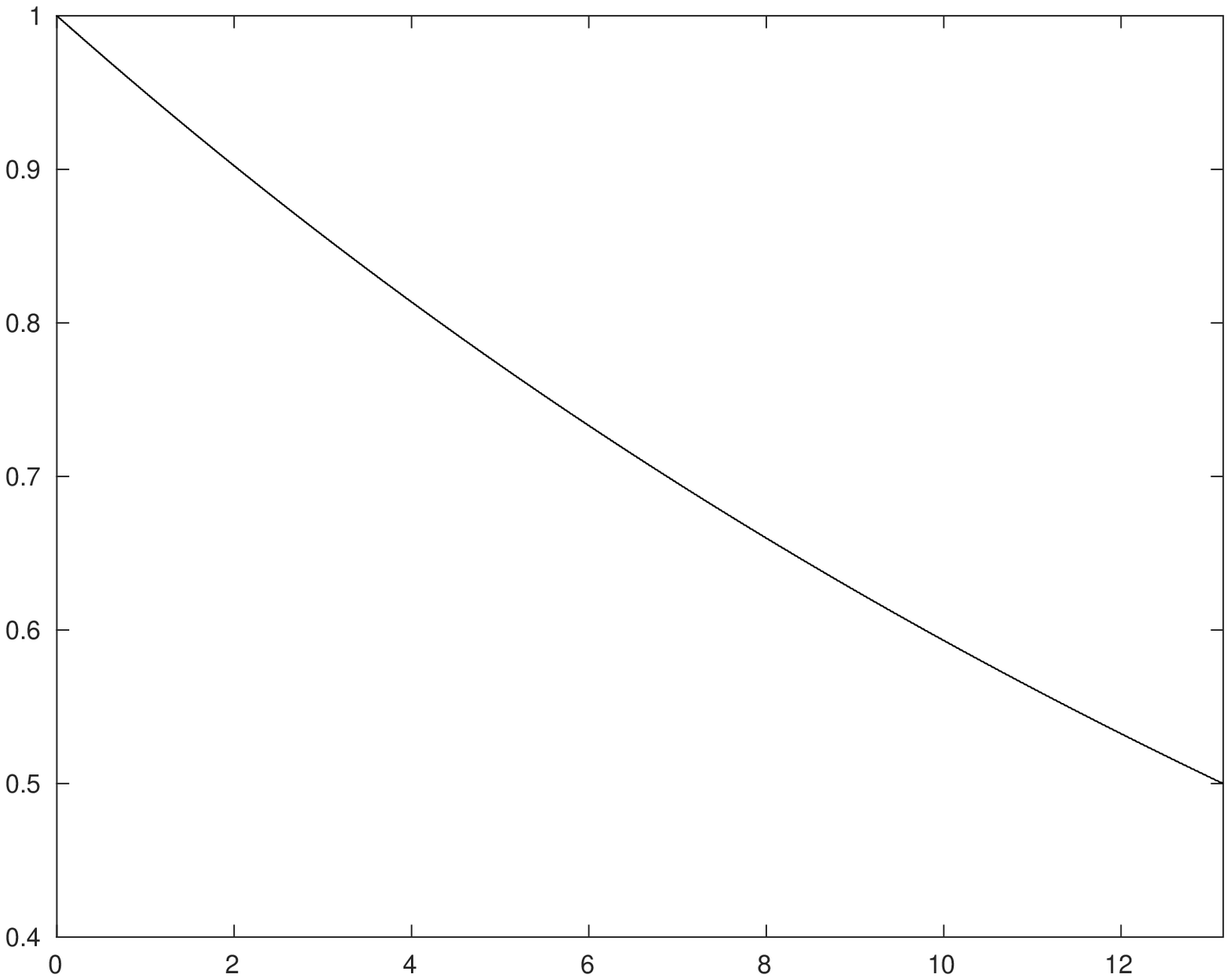}  
\hspace*{4mm}
\includegraphics[height=3.5cm]{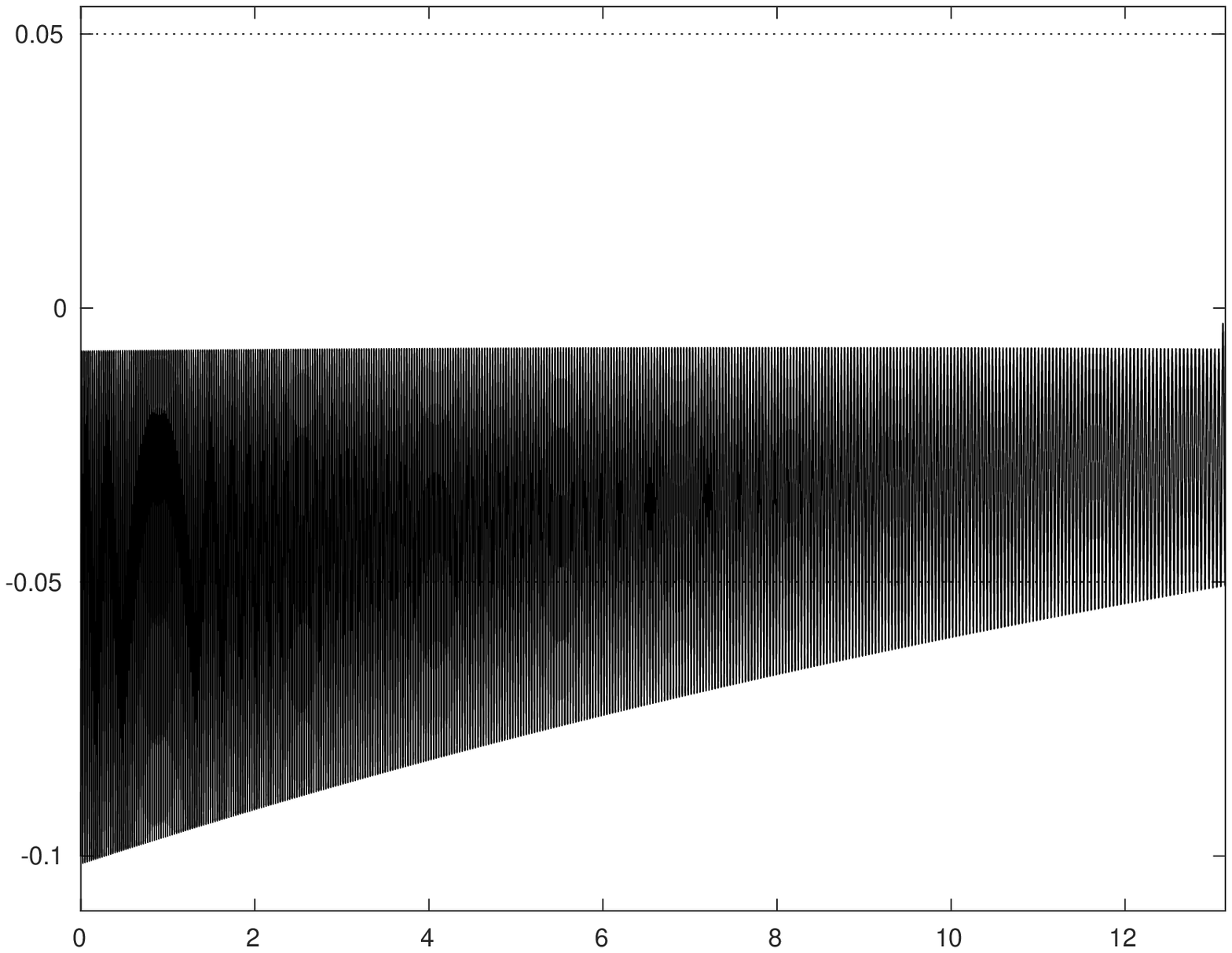}
\hspace*{4mm}
\includegraphics[height=3.5cm]{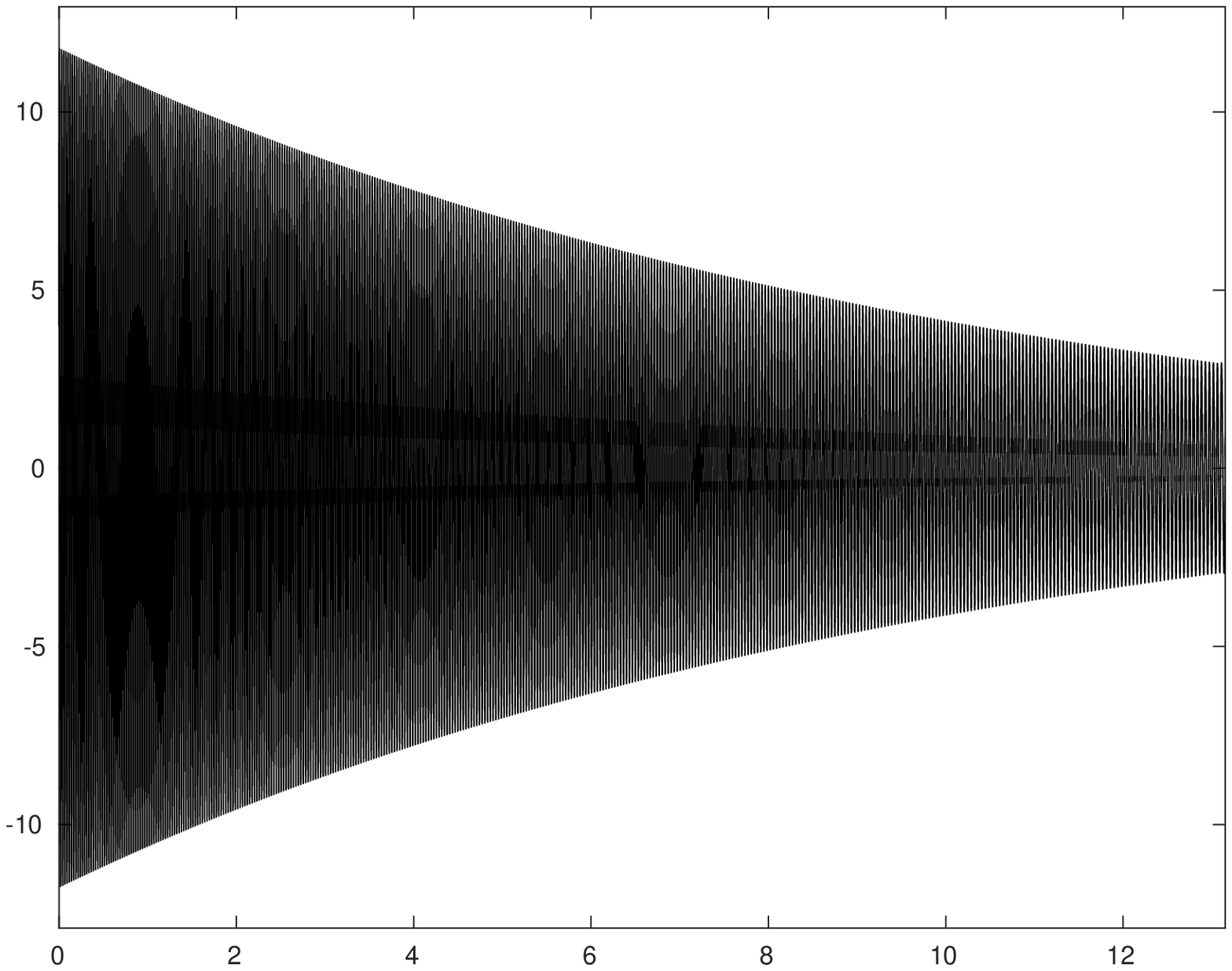}\\*[1.5ex]
\includegraphics[height=3.5cm]{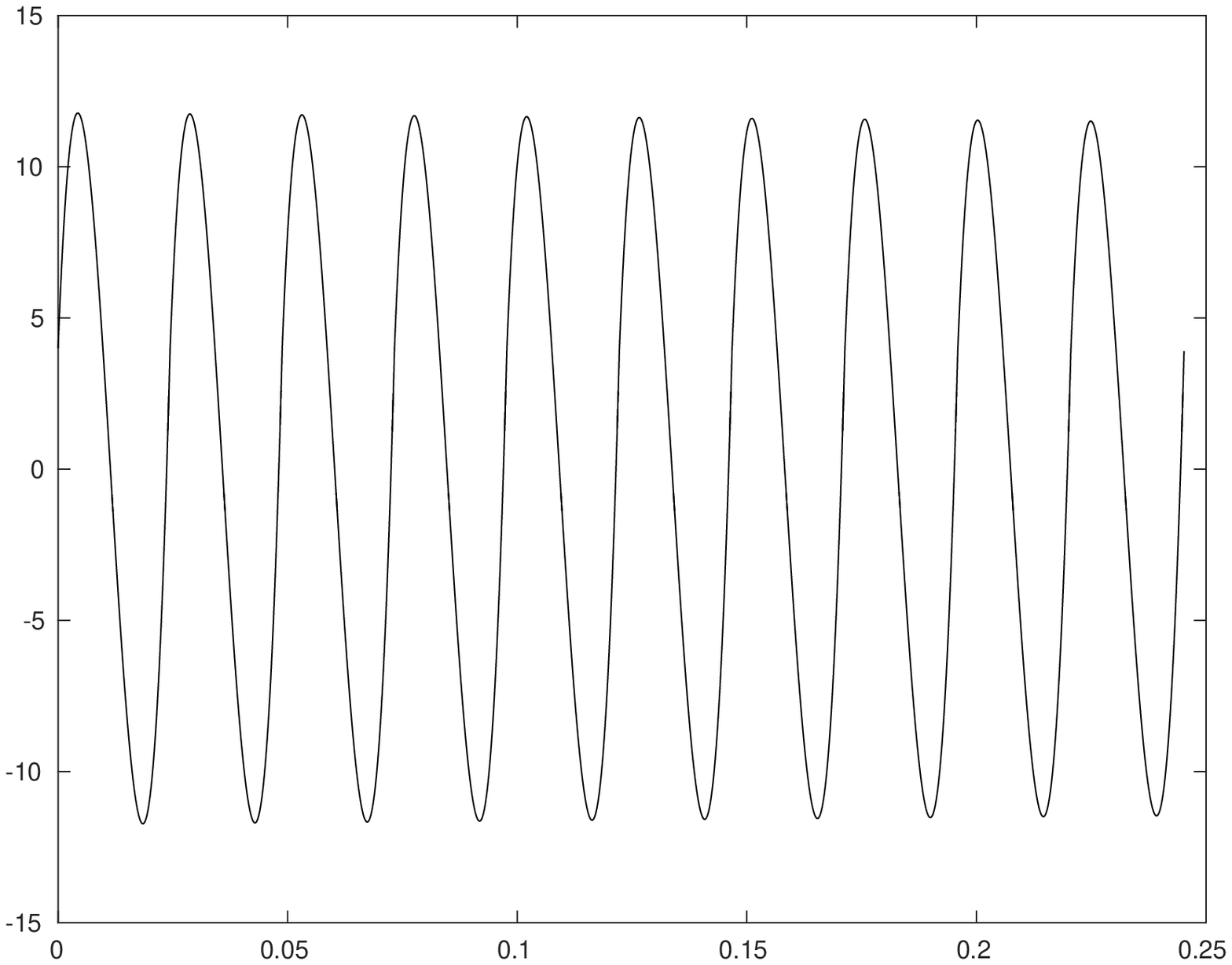}  
\hspace*{4mm}
\includegraphics[height=3.5cm]{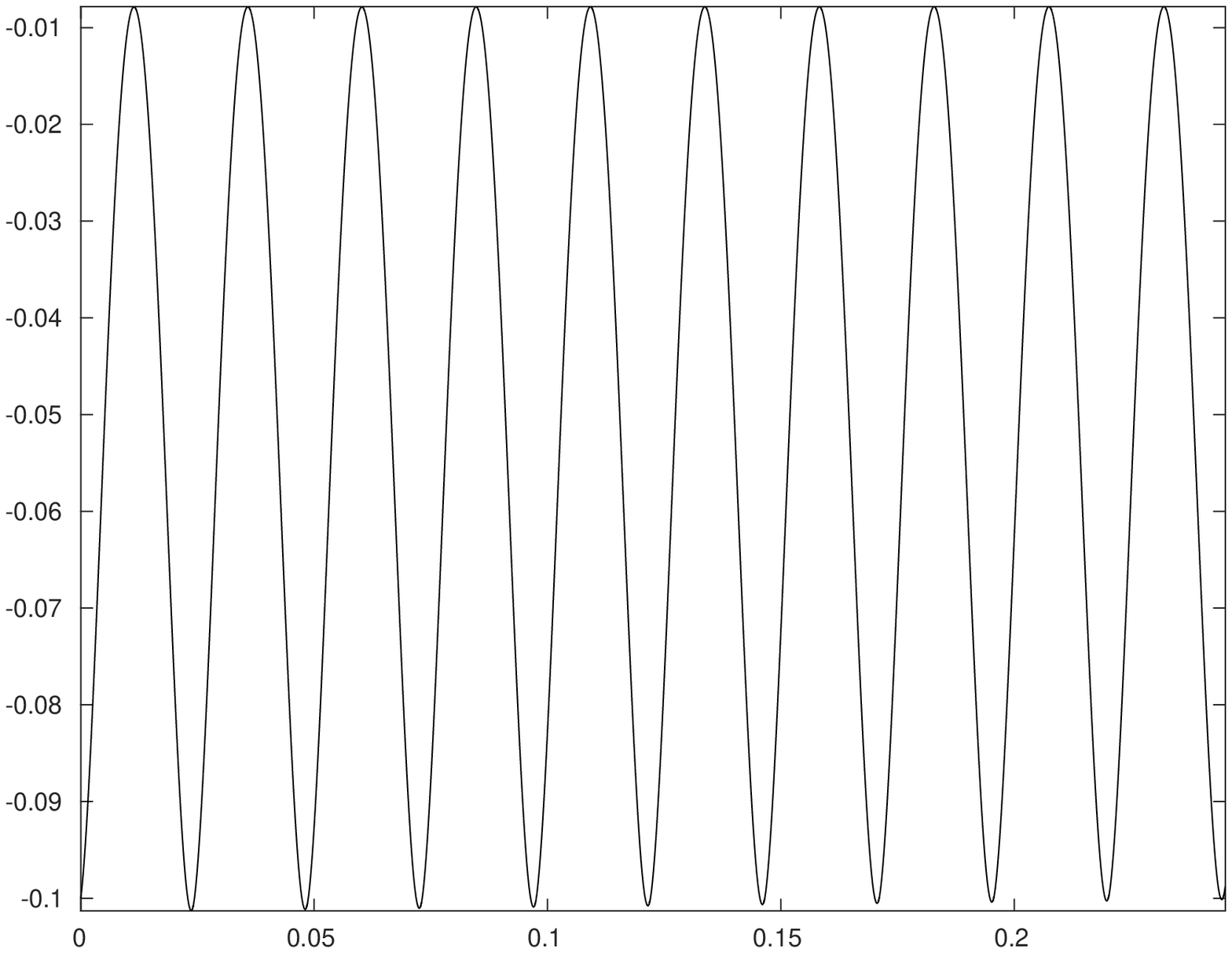}
\hspace*{4mm}
\includegraphics[height=3.5cm]{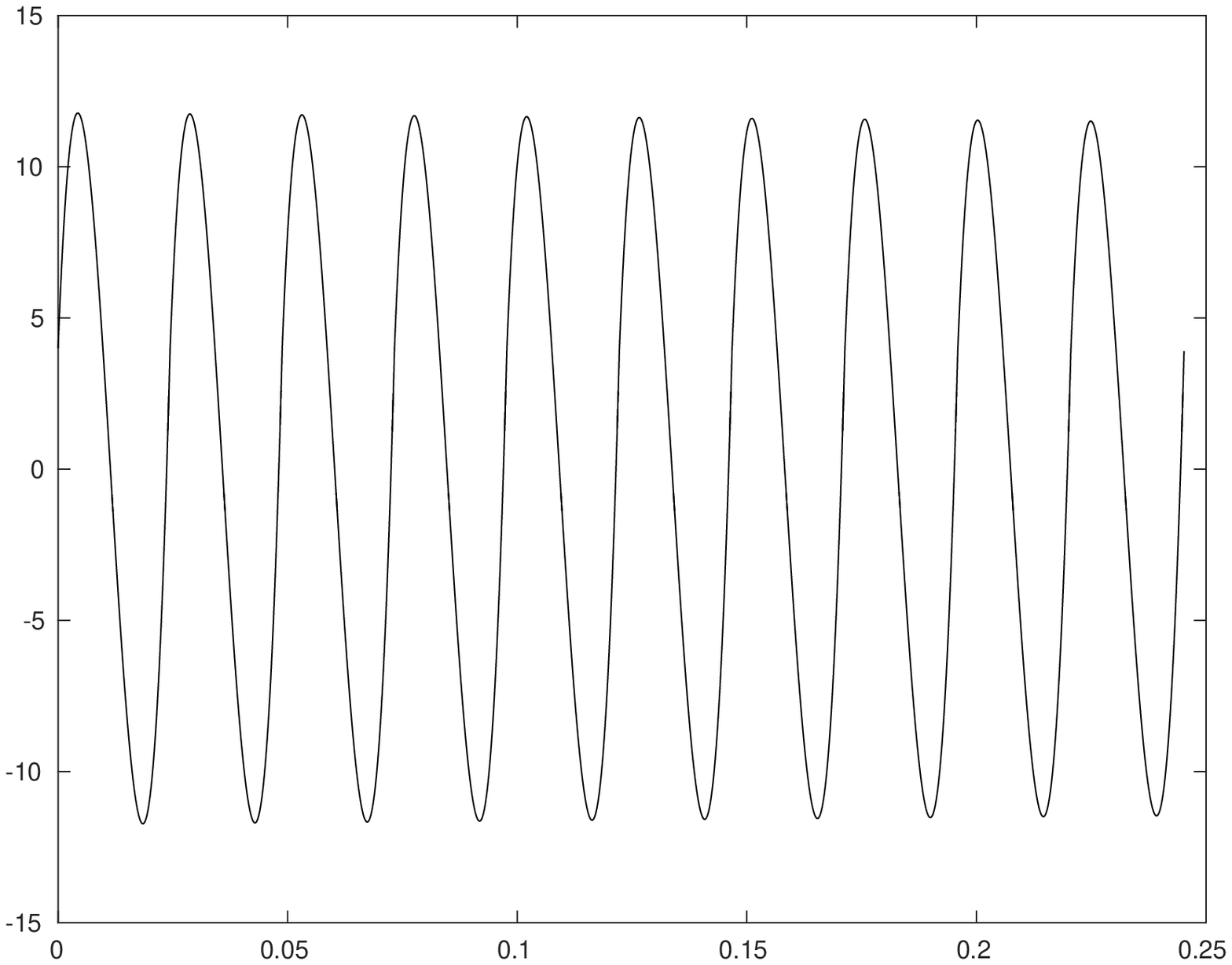}
\end{center}
\vspace*{-5mm}
\caption{\label{exalpha}$f^{\calM.\alpha}_\epsilon(x)$ (left) and its first
  (center) and second (right) derivatives as a function of $x$ for
  $\alpha = \half$ and $\epsilon = 5.10^{-2}$
  (top: $x \in [0,x_{k_{\epsilon,\alpha}}]$; bottom: $x \in  [0,x_{10}]$).
  Horizontal dotted lines indicate values of $-\epsilon$  
  and $\epsilon$ in the central top graph.}
\end{figure}

The top three graphs of Figure~\ref{exalpha} \vpageref{exalpha} illustrate the global behaviour
of the resulting function $f^{\calM.\alpha}_\epsilon(x)$ and of its first and second
derivatives for $x \in [0,x_{k_{\epsilon,\alpha}}]$, while the bottom ones show more detail of the first
10 iterations. The figure is constructed using $\epsilon = 5.10^{-2}$ and
$\alpha = \half$, which then yields that $k_{\epsilon,\alpha} = 148$. In
addition, we set $\lambda_k = \tenth |g_k|^{\frac{\alpha}{1+\alpha}}$ for
$k=0, \ldots, k_{\epsilon,\alpha}$. The nonconvexity of
$f^{\calM.\alpha}_\epsilon(x)$ is clear from the bottom graphs.

\llem{fAalpha-l}{The function $f^{\calM.\alpha}_\epsilon$ defined above on the
  interval $[0,x_{k_{\epsilon,\alpha}}]$ can be extended to a function from
  $\Re$ to $\Re$ satifying A.$\alpha$ and whose range is bounded independently
  of $\alpha$ and $\epsilon$.
}

\proof{
We start by showing that, on
\[
  [0,x_{k_{\epsilon,\alpha}}] = \bigcup_{k\in \calK} [x_k,x_k+s_k],
  \]
$f^{\calM.\alpha}_\epsilon$  is bounded in absolute value independently of
$\epsilon$ and $\alpha$, twice continuously differentiable with Lipschitz
continuous gradient and $\alpha$-H\"{o}lder continous Hessian. Recall first
\req{fM1} provide that $f^{\calM.\alpha}_\epsilon$ is twice continuously
differentiable by construction on $[0,x_{k_{\epsilon,\alpha}}]$. It thus
remains to investigate the gradient's Lipschitz continuity and Hessian's
$\alpha-$H\"older continuity, as well as whether $|f^{\calM.\alpha}_\epsilon(x)|$
is bounded on this interval.

Defining now
\beqn{ex-piphi-def}
\pi_k \eqdef \frac{\theta_k}{2}\frac{2f_k-1}{f_k}\in [0,\half \theta_k]
\tim{ and }
\phi(\theta)\eqdef 2 - \frac{1}{\theta} \in [2-\frac{\kappa_\lambda}{1-\kappa_{rg}},1+\kappa_{rg}]
\eeqn
(where we used \req{ex-fkbounds} and \req{ex-skdef}),
we obtain from \req{ex-nespalphadef}, \req{ex-fkgkHk}, \req{ex-skdef} and
\req{ex-expck}, that, for  $k\in\calK$,
\beqn{ex-ckbounds}
\begin{array}{llll}
  |c_{3,k}|s_k^2
  & \!= \epsilon f_k \left( 20 - \bigfrac{10}{\theta_k}-2\theta_k\right)
  -\epsilon^{\frac{3+2\alpha}{1+\alpha}}(4 + \pi_k)
  & \leq \epsilon\left[10|\phi(\theta)| + 2\theta + \frac{9}{2}\epsilon^{\frac{2+\alpha}{1+\alpha}}\right]
  & = \calO(\epsilon),\\%< 40 \,\epsilon \,\kappa_\lambda,\\
  |c_{4,k}|s_k^3
  & \!= \epsilon f_k \left(\bigfrac{15}{\theta_k}-30 + \theta_k\right)
   + \epsilon^{\frac{3+2\alpha}{1+\alpha}}(7 + 2\pi_k)
  & \leq \epsilon\left[15|\phi(\theta)| + \theta + 8\epsilon^{\frac{2+\alpha}{1+\alpha}}\right]
  & = \calO(\epsilon),\\%< 60 \, \epsilon \, \kappa_\lambda,\\
  |c_{5,k}|s_k^4
  & \!= \epsilon f_k \left(12 - \bigfrac{6}{\theta_k}\right)
  -\epsilon^{\frac{3+2\alpha}{1+\alpha}}(3 + \pi_k)
  & \leq \epsilon\left[6|\phi(\theta)| + \frac{7}{2}\epsilon^{\frac{2+\alpha}{1+\alpha}}\right]
  & = \calO(\epsilon),%< 24 \, \epsilon \, \kappa_\lambda.\\
\end{array}
\eeqn
where we also used $\epsilon \leq 1$ and \req{ex-fkbounds}. To show that the
Hessian of $f^{\calM.\alpha}_\epsilon$ is globally $\alpha-$H\"older continuous on
$[0,x_{k_{\epsilon,\alpha}}]$, we need to verify that \req{LipsHNEW} holds for
all $x,y$ in this interval. From \req{fM1}, this is implied by
\beqn{Hcont1}
|p^{'''}(s)|\leq c |s|^{-1+\alpha},\tim{for all $s\in [0,s_k]$ and $k \in \calK$,}
\eeqn
for some $c>0$ independent of $\epsilon$, $s$ and $k$. We have from the
expression of $p_k$ and $s\in [0,s_k]$ that 
\beqn{Hcont2}
\begin{array}{lcl}
|p_k^{'''}(s)|\cdot |s|^{1-\alpha}
& \leq & (6|c_{3,k}|+24|c_{4,k}|s_k+60 |c_{5,k}|s_k^2) s_k^{1-\alpha}\\*[1ex]
&   =  & (6|c_{3,k}|s_k^2+24|c_{4,k}|s_k^3+60 |c_{5,k}|s_k^4) s_k^{-(1+\alpha)}.
\end{array}
\eeqn
The boundedness of this last right-hand side on $[0,x_{k_{\epsilon,\alpha}}]$
, and thus the $\alpha$-H\"older continuity of the Hessian of $f^M$, then
follow from \req{ex-ckbounds}, \req{ex-skdef} and \req{ex-fkbounds}. 

Similarly, to show that the gradient of $f^M$ is globally Lipschitz continuous
in $[0,x_{k_{\epsilon,\alpha}}]$ is equivalent to proving that $p_k^{''}(s)$
is uniformly bounded above on the interval $[0,s_k]$ for $k\in \calK$. Since
$s_k>0$, we have 
\beqn{bounded-Hessian}
\begin{array}{lcl}
|p_k^{''}(s)|
& \leq& 2|c_{2,k}|+6|c_{3,k}|s_k+12|c_{4,k}|s_k^2+20|c_{5,k}|s_k^3\\*[1ex]
& = & 2|c_{2,k}| + (6|c_{3,k}|s_k^2+12|c_{4,k}|s_k^3+20|c_{5,k}|s_k^4)s_k^{-1}.
\end{array}
\eeqn
Then the third part of \req{ex-fkgkHk} and the bounds $\epsilon \leq 1$,
\req{ex-ckbounds}, \req{ex-c0c1}, \req{ex-skdef} and \req{ex-fkbounds} again
imply the boundedness of the last right-hand side on
$[0,x_{k_{\epsilon,\alpha}}]$, as requested. Finally, the fact that
$|f^{\calM.\alpha}_\epsilon|$ is bounded on $[0,x_{k_{\epsilon,\alpha}}]$ results
from the observation that, on the interval $[0,s_k]$ with $k\in \calK$,
\[
|p_k(s)| \leq f_k + |g_k||s_k| + \half |H_k|\,|s_k|^2
              + (|c_{3,k}|s_k^2+|c_{4,k}|s_k^3+|c_{5,k}|s_k^4)s_k
\]
from which a finite bound $a$ independent from $\alpha$ and $\epsilon$ again
follows from $\epsilon \leq 1$, \req{ex-fkgkHk}, \req{fM1}, \req{ex-ckbounds},
\req{ex-c0c1}, \req{ex-skdef} and \req{ex-fkbounds}. We have thus proved that
$f^{\calM.\alpha}_\epsilon$ satisfies the desired properties on
$[0,x_{k_{\epsilon,\alpha}}]$.

We may then smoothly prolongate $f^{\calM.\alpha}_\epsilon$ for $x\in \Re$, for
instance by defining two additional interpolation intervals
$[x_{-1},x_0]=[-1,0]$ and  $[x_{k_{\epsilon,\alpha}},x_{k_{\epsilon,\alpha}}+1]$
with end conditions
\[
f_{-1} = 1,
\ms
f_{k_{\epsilon,\alpha}+1} =
f_{k_{\epsilon,\alpha}}
\tim{and}
g_{-1}=H_{-1}=g_{k_{\epsilon,\alpha}+1} = H_{k_{\epsilon,\alpha}+1} = 0,
\]
and setting
\[
f^{\calM.\alpha}_\epsilon(x) = \left\{  \begin{array}{ll}
1 & \tim{for } x \leq -1,\\
p_k(x-x_k)+f_{k+1} & \tim{for } x\in[x_k,x_{k+1}]
                                 \tim{and} k\in \iibe{-1}{k_{\epsilon,\alpha}},\\
f^{\calM.\alpha}_\epsilon(x_{k_{\epsilon,\alpha}}) & \tim{for } x \geq x_{k_{\epsilon,\alpha}}+1,\\
\end{array}\right.
\]
which subsumes \req{fM1}. Using arguments similar to those used above, it is
easy to verify from \req{ex-c0c1}, \req{ex-expck}  and $s_{-1} =
s_{k_{\epsilon,\alpha}}=1$ that all desired properties are maintained.
}

\noindent
We formulate the results of this development in the following theorem.

\lthm{slowMa}{ For every $\epsilon \in (0,1)$, every $\alpha \in [0,1]$ and
  every method in $\calM.\alpha$, a function $f^{\calM.\alpha}_\epsilon$ satisfying
  A.$\alpha$ with values in a bounded interval independent of $\epsilon$ and
  $\alpha$ can be constructed, such, when applied to
  $f^{\calM.\alpha}_\epsilon$, the considered method terminates 
  exactly at iteration
  \[
  k_{\epsilon,\alpha}=\left\lceil\epsilon^{-\frac{2+\alpha}{1+\alpha}}\right\rceil.
  \]
  with the first iterate $x_{k_{\epsilon,\alpha}}$ such that
  $\|\nabla_x f^{\calM.\alpha}_\epsilon(x_{k_{\epsilon,\alpha}})\| \leq \epsilon$.
}

Note that the prolongation of $f^{\calM.\alpha}_\epsilon(x)$ to $x\geq0$ suggested
as an example in the proof of Lemma~\ref{fAalpha-l} admits an isolated finite
global minimizer. Indeed, since the $g_{k_{\epsilon,\alpha}}<0$, there must be
a value lower than $f(x_{k_{\epsilon,\alpha}})$ in
$(x_{k_{\epsilon,\alpha}},x_{k_{\epsilon,\alpha}}+1)$, and thus the global
minimizer must lie in one of the constructed sub-intervals in
$(-1,x_{k_{\epsilon,\alpha}+1})$; since $f^{\calM.\alpha}_\epsilon(x)$ is quintic
(and not constant) in each of these, the global minimizer must therefore be
isolated.

\subsection{The inexact Newton's method}\label{slowN-s}

It is interesting that the technique developed in the previous subsection can
also be used to derive an $\calO\left(\epsilon^{-2}\right)$ lower bound on
worst-case evaluation complexity for an inexact Newton's method applied to a function
having Lipschitz continuous Hessians on the path of iterates.  This is stronger
than using Theorem~\ref{slowMa} above for $\alpha = 1$, as it would result in
a weaker $\calO\left(\epsilon^{-3/2}\right)$ lower bound, or for $\alpha = 0$ as
it would then only guarantee bounded Hessians. In the spirit of \cite{CartGoulToin10a}, this
new function is constructed by extending  to $\Re^2$ the unidimensional 
$f^{\calM.0}_\epsilon(x)$ obtained in the previous section for the specific
choice $M_k=0$, which then ensures that $\theta_k\in
[1-\kappa_{rg},1+\kappa_{rg}]$ for all $k$ (see \req{ex-lambdabound} and
\req{ex-skdef}). The proposed extension is of the form
\beqn{h-def}
h^N_\epsilon(x,y) \eqdef f^{\calM.0}_\epsilon(x) + u_\epsilon(y),
\eeqn
where we still have to specify the univariate function $u_\epsilon$ such
that Newton's method applied to $u_\epsilon$ converges with large steps. In
order to define it, we start by redefining
\[
k_\epsilon = k_{\epsilon,0} = \lceil \epsilon^{-2} \rceil
\tim{ and }
\calK = \iibe{0}{k_\epsilon}.
\]
Then we set, for $k\in \calK$,
\beqn{new-fkgkHk-y}
u_k = 1 - \half k \epsilon^2,
\ms
g_k^u = - 2\epsilon^2 u_k,
\ms
H_k^u = 2|g_k^u|u_k > 0,
\eeqn
and
\beqn{sku-def}
s_k^u = \frac{\nu_k}{2u_k}
\tim{with} \nu_k \in \left[1-\kappa_{rg},1+\kappa_{rg}\right]
\tim{and}
u_k \in [\half,1],
\eeqn
this definition allowing for
\[
H_k^u s_k^u = -g_k^u + r_k^u
\tim{ with }
|r_k^u| \leq \kappa_{rg}|g^u_k|.
\]
(Remember that $M_k=0$ because we are considering Newton's method.)
Note that sufficient decrease is obtained in manner similar to \req{ex-d1}-\req{ex-decrease},
because of \req{new-fkgkHk-y}, \req{sku-def} and $\lambda_k=0$, yielding that
$u_k-u_{k+1} \geq -(g_k^us_k^u + \half H_k^u (s_k^u)^2) /(1+\kappa_{rg})$.
Setting now $y_0= 0$ and $y_{k+1}= y_k + s_k^u$ for $k \in \ii{k_\epsilon}$, we
may then, as in Section~\ref{slowMa-s}, define
\beqn{fM2}
u_\epsilon(y)=  p^u_k(y-y_k)+u_{k+1}
\tim{for $y\in[y_k,y_{k+1}]$ and $k = 0,\ldots, k_\epsilon-1$,}
\eeqn
where $p^u_k$ is a fifth degree polynomial interpolating
the values and derivatives given by \req{new-fkgkHk-y} on the
interval $[0,s_k^u]$. We then obtain the following result.

\lthm{slowN}{ For every $\epsilon \in (0,1)$,
  there exists a function $h^N_\epsilon$ with
  Lipschitz continuous gradient and Lipschitz continuous Hessian along the
  path of iterates 
%  \[
  $\cup_{k = 0}^{k_\epsilon-1} [x_j,x_{j+1}]$,
%  \]
  and with values in a bounded interval independent of $\epsilon$, such that,
  when applied to $h^N_\epsilon$, Newton's terminates exactly at iteration
  \[
  k_\epsilon=\left\lceil\epsilon^{-2}\right\rceil
  \]
  with the first iterate $x_{k_\epsilon}$ such that
  $\|\nabla_x f^{\calM.\alpha}_\epsilon(x_{k_\epsilon})\| \leq \epsilon\sqrt{1+\epsilon^2}$.
}

\proof{
One easily verifies from \req{new-fkgkHk-y}, \req{sku-def} and 
\req{ex-expck} that the interpolation coefficients, now denoted by
$|d_{i,k}|$, are bounded for all $k \in \iibe{0}{k_\epsilon-1}$ and
$i \in \iibe{0}{5}$. This observation and \req{sku-def} in turn guarantee
that $u_\epsilon$ and all its derivatives (including the third)
remain bounded on each interval $[0,s_k^u]$ by constants independent of
$\epsilon$. As in Lemma~\ref{fAalpha-l}, we next extend
$u_\epsilon$ to the whole of $\Re$ while preserving this property. We then
construct $h^N$ using \req{h-def}. From the properties of $f^{\calM.0}_\epsilon$ and
$u_\epsilon$, we deduce that
$h^N_\epsilon$ is twice continuously differentiable and has
a range bounded independently of $\epsilon$.  Moreover, it
satisfies A.0.  When applied on $h^N_\epsilon(x,y)$, Newton's generates the
iterates $(x_k,y_k)$ and its gradient at the $k_\epsilon$-th iterate is
$( \epsilon, \epsilon^2)$ so that
$\|\nabla h^N(x_{k_\epsilon},y_{k_\epsilon})\| = \epsilon\sqrt{1+\epsilon^2}$,
prompting termination. Before that, the algorithm generates the steps
$(s_k,s_k^u)$, where, because both $f_k$ and $u_k$ belong to
$[\half,1]$ and because of \req{ex-skdef} with $\alpha= 0$, 
\beqn{skubounds}
s_k \in [\epsilon(1-\kappa_{rg}),2\epsilon(1+\kappa_{rg})]
\tim{ and }
s_k^u \in [1-\kappa_{rg},2(1+\kappa_{rg})].
\eeqn
Thus the absolute value
of the third derivative of $h^N_\epsilon(x,y)$ is given, for
$(x,y)$ in the $k$-th segment of the path of iterates, by
\beqn{3rdder}
\begin{array}{lcl}
\lefteqn{\frac{1}{\|(s_k,s_k^u)\|}\Big| p_k^{'''}(x-x_k) s_k^3+ (p_k^u)^{'''}(y-y_k)(s_k^u)^3\Big|}&&\\*[1.2ex]
& \ms\ms \leq & \bigfrac{1}{1-\kappa_{rg}}\,\Big[|p_k^{'''}(x-x_k)|s_k^3 + |(p_k^u)^{'''}(y-y_k)|(s_k^u)^3\Big]\\*[1.2ex]
& \ms\ms   =  & \bigfrac{1}{1-\kappa_{rg}}\,\Big[\,\Big(6|c_{3,k}|+24|c_{4,k}|s_k+60|c_{5,k}|s_k^2\Big) s_k^3 \\*[1.2ex]
&             & \ms\ms\ms\;\; + \Big(6|d_{3,k}|+24|d_{4,k}|s_k^u+60 |d_{5,k}|(s_k^u)^2\Big)(s_k^u)^3\Big]\\*[1.2ex]
& \ms\ms   =  & \bigfrac{1}{1-\kappa_{rg}}\,\Big[\,\Big(6|c_{3,k}|s_k^2+24|c_{4,k}|s_k^3+60|c_{5,k}|s_k^4\Big) s_k \\*[1.2ex]
&             & \ms\ms\ms\;\; + 6|d_{3,k}|(s_k^u)^3+24|d_{4,k}|(s_k^u)^4+60 |d_{5,k}|(s_k^u)^5\Big],
\end{array}
\eeqn
where we used the fact that
$\|(s_k,s_k^u)\|\geq \|s_k^u\|.$ and \req{skubounds}.
But, in view of
\req{ex-ckbounds}, \req{ex-piphi-def} with $\theta_k\in
    [1-\kappa_{rg},1+\kappa_{rg}]$, \req{skubounds}, 
$\epsilon\leq 1$ and the boundedness of the $d_{i,k}$, the last right-hand 
side of \req{3rdder} is bounded by a constant independent of $\epsilon$.  Thus
the third derivative of $h^N_\epsilon(x,y)$ is bounded on every segment by the
same constant, and, as a consequence, the Hessian of $h^N_\epsilon(x,y)$ is
Lipschitz continuous of each segment, as desired.
}

\noindent
Note that the same result also holds for any method in $\calM.0$ with
$M_k$ small enough to guarantee that $s_k$ is bounded away from zero for all $k$.

\numsection{Complexity and optimality for methods in $\calM.\alpha$}

We now consider the consequences of the examples derived in Section~3 on the
evaluation complexity analysis of the various methods identified in Section~2 as
belonging to $\calM.\alpha$.

\subsection{Newton's method.}

First note that the third part of \req{ex-fkgkHk} ensures that $H_k>0$ so that
the Newton iteration is well-defined for the choice \req{Newton}.  This choice
corresponds to setting $\theta_k = 1$ for all $k \geq 0$ in the example of
Section~3. So we first conclude from Theorem~\ref{slowMa} that Newton's method
may require $\epsilon^{-(2+\alpha)/(1+\alpha)}$ evaluations when applied on
the resulting objective function $f^{\calM.\alpha}_\epsilon$ satisfying
A.$\alpha$ to generate $|g_k|\leq \epsilon$. However, Theorem~\ref{slowN}
provides the stronger result that it may in fact require $\epsilon^{-2}$
evaluations (as a method in $\calM.0$) for nearly the same task (we traded
Lipschitz continuity of the Hessian on the whole space for that along the path
of iterates). As a consequence we obtain that \emph{Newton's method is not
optimal in $\calM.\alpha$ as far as worst-case evaluation complexity is
concerned}.

The present results also improves on the similar bound given in \cite{CartGoulToin11c},
in that the objective function on Sections~\ref{slowMa-s} and \ref{slowN-s}
ensure the existence of a lower bound $f_{\rm low}$ on
$f^{\calM.\alpha}_\epsilon(x)$ such that $f^{\calM.\alpha}_\epsilon(x_0) - f_{\rm
  low}$ is bounded, while the latter difference is unbounded in \cite{CartGoulToin11c}
(for $\alpha \in \{0,1\}$) as the number of iterations approaches
$\epsilon^{-2}$. We will return to the significance of this observation when
discussing regularization methods.

Since the steepest-descent method is known to have a worst-case evaluation
complexity of $\calO\left(\epsilon^{-2}\right)$ when applied on functions having
Lipschitz continuous gradients \cite[p.~29]{Nest04} , Theorem~\ref{slowN}
shows that Newton's method may, in the worst case, converge as slowly as
steepest descent in the worst case. Moreover, we show in Appendix A1 that the
quoted worst-case evaluation complexity bound for steepest descent is sharp,
which means that steepest-descent and Newton's method are undistinguishable
from the point of view of worst-case complexity orders.

Note also that if the Hessian of the objective is unbounded, and hence, we are
outside of the class A.$0$, the worst-case evaluation complexity of Newton's
method worsens, and in fact, it may be arbitrarily bad \cite{CartGoulToin10a}.

\subsection{Cubic and other regularizations.}

Recalling our discussion of the $(2+\alpha)$-regularization method in
Section~\ref{NtrARC}, we first note, in the example of Section~\ref{slowMa-s},
that, because of \req{skdef} and \req{LAMBDAprop},
$s_k$ is a minimizer of the model \req{model} with $\beta_k= \lambda_k$ at
iteration $k$, in that 
\beqn{ex-success}
m_k(s_k) = f^{\calM.\alpha}_\epsilon(x_k+s_k)= f_{k+1}
\eeqn
for $k \in \calK$.  Thus every iteration is successful as the
objective function decrease exactly matches decrease in the model. Hence the
choice $\sigma_k=\sigma>0$ for all $k$ is allowed by the method, and thus
$\lambda_k= \sigma\|s_k\|^{2+\alpha}$ satisfies \req{LAMBDAprop} and
\req{lambdacondAL}. Theorem~\ref{slowMa} then shows that 
this method may require at least $\epsilon^{-(2+\alpha)/(1+\alpha)}$
iterations to generate an iterate with $|g_k| \leq \epsilon$.  This is
important as the \emph{upper} bound on this number of iterations was
proved\footnote{As a matter of fact, \cite{CartGoulToin11d}   contains a
detailed proof of the result for $\alpha=1$, as well as the statement 
that it generalizes for  $\alpha \in (0,1]$.  Because of the central role
of this result in the present paper, a more detailed proof of the worst-case
evaluation complexity bound for $\alpha \in (0,1]$ in provided as
Appendix~A2.} in \cite{CartGoulToin11d} to be
\beqn{upper-reg}
O\Big([ f(x_0)-f_{\rm low})] \, \epsilon^{-\frac{2+\alpha}{1+\alpha}} \Big)
\eeqn
where $f_{\rm low}$ is any lower bound of $f(x)$.
Since we have that $f(x_0)-f_{\rm low}$ is a fixed number independent of
$\epsilon$ for the example of Section~\ref{slowMa-s}, this shows that the ratio
\beqn{rho-comp}
\rho_{\rm comp}
\eqdef \frac{\tim{upper bound on the worst-case evaluation complexity}}
            {\tim{lower bound on the worst-case evaluation complexity}}
\eeqn
for the $(2+\alpha)$-regularization method is bounded independently of
$\epsilon$ and $\alpha$. Given that \req{upper-reg} involves an unspecified
constant, this is the best that can be obtained as far as the
order in $\epsilon$ is concerned, and yields the following important result
on worst-case evaluation complexity.

\lthm{reg-sharp-th}{When applied to a function satisfying A.$\alpha$, the 
  $(2+\alpha)$-regularization method may require at most \req{upper-reg} 
  function and derivatives evaluations.  Moreover this bound is sharp (in the
  sense that $\rho_{\rm comp}$ is bounded independently of $\epsilon$ and
  $\alpha$) and the $(2+\alpha)$-regularization method is optimal in $\calM.\alpha$.
}

\proof{The optimality of the $(2+\alpha)$-regularization method within
  $\calM.\alpha$ results from the observation that the example of Section~3
  implies that no method in $\calM.\alpha$ can have a worst-case evaluation
  complexity of a better order.}

\noindent
In particular, the cubic regularization method is optimal for smooth
optimization problems with Lipschitz continuous second derivatives.  As we
have seen above, this is in contrast with Newton's method.

Note that Theorem~\ref{reg-sharp-th} as stated does \emph{not} result from the
statement in \cite{CartGoulToin11c} that the bound
\req{upper-reg} is ``essentially sharp''. Indeed this latter statement
expresses the fact that, for any $\tau>0$, there exists a function independent
of $\epsilon$, on which the relevant method may need at least
$\epsilon^{-3/2+\tau}$ evaluations to terminate with $|g_k|\leq
\epsilon$. But, for any fixed $\epsilon$, the value of $f(x_0)-f_{\rm low}$
tends to infinity when, in the example of that paper, the number of iterations
to termination approaches $\epsilon^{-3/2}$ as $\tau$ goes to zero. As a
consequence, the numerator of the ratio \req{rho-comp}, that is
\req{upper-reg}, and $\rho_{\rm comp}$ itself are unbounded for that
example. Theorem~\ref{reg-sharp-th} thus brings a formal improvement on the
conclusions of \cite{CartGoulToin11c}.

\subsection{Goldfeld-Quandt-Trotter}

Recalling \req{gqtCONSTR}, we can set $\omega_k=\omega$ in the algorithm as
every iteration is successful due to \req{ex-success} which, with
\req{ex-fkgkHk} and $f_k \in [\half, 1]$ gives that $\lambda_k +
\lambda_{\min}(H_k)\leq \omega |g_k|^{\frac{\alpha}{1+\alpha}}$, which is in
agreement with \req{sklbdUPPER} and \req{lambdacondAL}.  Thus the lower bound
of $\epsilon^{-(2+\alpha)/(1+\alpha)}$ iterations for termination also applies
to this method.

An upper bound on the worst-case evaluation complexity for the GQT method can
be obtained by the following argument. We first note that, similarly to
regularization methods, we can bound the total number of unsuccessful
iterations as a constant multiple of the successful ones, provided $\omega_k$
is chosen such that \req{Rmin} holds. Moreover, since $f$ satisfies
A.$\alpha$, its Hessian is bounded above by \req{Hbounded}.  In addition, we
have noted in Section~\ref{NtrARC} that $\|g_k\|$ is also bounded above.  In
view of \req{gqtCONSTR} and \req{Rbd}, this in turn implies that $\|H_k +
\lambda_kI\|$ is also bounded above.  Hence we obtain from \req{GQT-step} that
$\|s_k\| \geq \kappa_{GQT} \|g_k\| \geq \kappa_{GQT} \,\epsilon$ for some
$\kappa_{QGT}>0$, as along as termination has not occurred.  This last bound
and \req{QGTmodeldecrease} then give that GQT takes at most
$\calO\left((f(x_0)-f_{\rm low})\epsilon^{-\frac{\alpha}{1+\alpha}-2}\right)$
iterations, which is worse than \req{upper-reg} for $\alpha>0$.  Note that
this bound improves if only Newton steps are taken (i.e. $\lambda_k = 0$ is
chosen for all $k\geq 0$), to be of the order of \req{upper-reg}; however,
this cannot be assumed in the worst-case for nonconvex functions. In any case,
it implies that the GQT method is not optimal in $\calM.\alpha$.

\subsection{Trust-region methods}

Recall the choices \req{trust-region} we make in this case. If $\lambda_k=0$,
the trust-region constraint $\|s\|\leq \Delta_k$ is inactive at $s_k$, in
which case, $s_k$ is the Newton step. If we make precisely the choices we made
for Newton's method above, choosing $\Delta_0$ such that $\Delta_0>|s_0|$
implies that the Newton step will be taken in the first and in all subsequent
iterations since each iteration is successful and then $\Delta_k$ remains
unchanged or increases while the choice \req{ex-skdef} implies that $s_k$
decreases. Thus the trust-region approach, through the Newton step, has a
worst-case evaluation complexity when applied to $f^{\calM.\alpha}_\epsilon$ which
is at least that of the Newton's method, namely $\epsilon^{-2}$.

\comment{%%%%%%%%%%%%%%%%%%%%%%%%%%%%%%%%%%%%%%%%%%%%%%%%%%%%%%%%%%%%%%%%%%%%%%
{\bf Do we care for the next (somewhat confusing) paragraph???  The previous
  one is enough to show that TR's complexity is at least as bad as Newton's.}

By contrast when $\lambda_k> 0$ for all $k$,  $s_k=\Delta_k$. Using the
notation  in \cite[Algorithm 6.1.1]{ConnGoulToin00}, let $\eta$ in
\req{fmsuccess} be equal to $\eta_1$, which corresponds to successful but not
very successful steps $s_k$. This allows  the trust-region radius $\Delta_k$
to decrease slightly, namely $\Delta_{k+1}\in [\gamma_2\Delta_k,\Delta_k]$. We
let  
\[
\Delta_{k+1}=\gamma_k\Delta_k, 
\tim{where $\gamma_k=\left(\frac{k+1}{k+2}\right)^{q}$} 
\]
with $q$ defined in \req{q3}. If $\gamma_2$ is chosen such that $\gamma_2\leq
(1/2)^q$, then clearly the above updating rule implies that $\Delta_{k+1}\in
[\gamma_2\Delta_k,\Delta_k]$. Note that due to \req{skex}, this updating rule
is consistent with the trust-region constraint being active on each
iteration. Now we can choose $H_k<0$ so as to ensure \req{H}, noting that
despite not knowing the precise value of $\lambda_k$ for the trust-region
method, we know that the global solution of the trust-region subproblem  is
unique whenever $H_k+\lambda_k I$ is positive definite, which is clearly the
case here, due to \req{H}.
}%%%%%%%%%%%%%%%%%%%%%%%%%%%%%%%%%%%%%%%%%%%%%%%%%%%%%%%%%%%%%%%%%%%%%%%%%%%%%%

\subsection{Linesearch methods}

Because the examples of Sections~\ref{slowMa-s} and \ref{slowN-s} are valid
for $r_k=0$ which corresponds to $\mu_k=1$ for all $k$, and because this
stepsize is acceptable since $f(x_{k+1}) = m_k(s_k)$, we 
deduce that at least $\epsilon^{-\frac{2+\alpha}{1+\alpha}}$ iterations and
evaluations may be needed for the linesearch variants of any method in
$\calM.\alpha$ applied to a function satisfying A.$\alpha$, and that
$\epsilon^{-2}$ evaluations may be needed for the linesearch variant of
Newton's method applied on a function satisfying A.$0$.  Thus the conclusions
drawn regarding their (sub-)optimality in terms of worst-case evaluation
complexity are not affected by the use of a linesearch.

\numsection{The Curtis-Robinson-Samadi class}

We finally consider a class of methods recently introduced in
\cite{CurtRobiSama17c}, which we call the CRS class. 
This class depends on the parameters $0<\underline{\sigma}\leq\bar{\sigma}$,
$\eta \in (0,1)$ and two non-negative accuracy thresholds $\kappa_1$ and
$\kappa_2$. It is defined as follows. At the start, adaptive regularization
thresholds are set according to 
\beqn{CRS-sigma-bounds}
\sigma^L_0=0 \tim{and} \sigma^U_0= \bar{\sigma}.
\eeqn
Then for each iteration $k\geq 0$, a step $s_k$ from the current iterate $x_k$
and a regularization parameter $\lambda_k\geq 0$ are chosen to
satisfy\footnote{In \cite{CurtRobiSama17c}, further restrictions on the step
  are imposed in order to obtain global convergence under A.0 and bounded
  gradients, but are irrelevant for the worst-case complexity analysis 
  under A.$1$. We thus ignore them here, but note that this analysis also
  ensures global convergence to first-order stationary points.}
\beqn{CRS-inexact}
(H_k + \lambda_k I) s_k = -g_k + r_k,
\eeqn
\beqn{CRS-lambda-bounds}
\sigma_k^L \|s_k\| \leq \lambda_k \leq \sigma_k^U \|s_k\|,
\eeqn
\beqn{CRS-slope}
s_k^Tr_k \leq \half s_k^T(H_k+\lambda_k I)s_k + \half \kappa_1\|s_k\|^3,
\eeqn
and
\beqn{CRS-srule}
\|r_k\| \leq \lambda_k\|s_k\| + \kappa_2 \|s_k\|^2.
\eeqn
The step is then accepted, setting $x_{k+1}=x_k+s_k$, if
\beqn{CRS-rho}
\rho_{CRS} = \frac{f(x_k)-f(x_k+s_k)}{\|s_k\|^3} \geq \eta
\eeqn
or rejected otherwise.  In the first case, the regularization thresholds are
reset according to \req{CRS-sigma-bounds}.  If $s_k$ is rejected, $\sigma^L_k$
and $\sigma^U_k$ are updated by a simple mechanism (using
$\underline{\sigma}$) which is irrelevant for our purpose here.  The algorithm
is terminated as soon as an iterate is found such that $\|g_k\| \leq \epsilon$.

Observe that \req{CRS-inexact} corresponds to inexactly minimizing the
regularized model \req{model} and that \req{CRS-srule} is very similar
to the subproblem termination rule of \cite{BirgGardMartSantToin17}.

An upper bound of $\calO\left(\epsilon^{-3/2}\right)$ is proved in
\cite[Theorem~17]{CurtRobiSama17c} for the worst-case evaluation complexity of
the methods belonging to the CRS class.  It is stated in
\cite{CurtRobiSama17c} that both ARC 
\cite{Grie81,WeisDeufErdm07,NestPoly06,CartGoulToin11,CartGoulToin11d}
and TRACE \cite{CurtRobiSama17} belong to the class, although the details are
not given.

Clearly, the CRS class is close to $\calM.1$, but yet differs from it.  In
particular, no requirement is made that $H_k+\lambda_k I$ be positive
semi-definite but \req{CRS-slope} is required instead, there is no formal
need for the step to be bounded and \req{CRS-srule} combined with
\req{CRS-lambda-bounds} is slightly more permissive than the second part of
\req{skdef}. We now define CRS$_a$, a sub-class of the CRS class of methods,
as the set of CRS methods for which \req{CRS-srule} is
strengthened\footnote{Hence the subscript $a$, for ``accurate''.} to become 
\beqn{CRS-srule2}
\|r_k\| \leq \min\Big[\kappa_{rg}\|g_k\|,\lambda_k\|s_k\|+\kappa_2\|s_k\|^2 \Big]
\tim{ with } \kappa_{rg} < 1.
\eeqn
(in a manner reminiscent of the second part of \req{skdef}) and such that
\beqn{CRS-eta-cond}
2\eta(1+\kappa_{rg})^3 \leq 1
\eeqn
(a mild technical condition\footnote{Due to the lack of scaling invariance of
\req{CRS-rho}, at variance with \req{rho}.} whose need will become apparent
below). We claim that, for any choice of method in the CRS$_a$ class and
termination threshold $\epsilon$, we can construct a function satisfying A.1
such that the considered CRS$_a$ method terminates in exactly
$\left\lceil\epsilon^{-3/2}\right\rceil$ iterations and evaluations.
This achieved simply by showing that the generated sequences of iterates,
function, gradient and Hessian values belong to those detailed in
the example of Section~\ref{slowMa-s}.

We now apply a method of the CRS$_a$ class for a given $\epsilon >0$, and
first consider an iterate $x_k$ with associated values $f_k$, $g_k$
and $H_k$ given by \req{ex-fkgkHk} for $\alpha=1$, that is
\beqn{CRS-fkgkHk}
f_0 = 1,
\ms
f_k = f_0 - \half k \epsilon^{3/2},
\ms
g_k = -2 \epsilon f_k
\tim{and}
H_k = 4 \epsilon^{1/2}f_k^2;
\eeqn
Suppose that
\beqn{CRS-sigmakb}
\sigma_k^L=0 \tim{ and } \sigma_k^U = \bar{\sigma}
\eeqn
(as is the case by definition for $k=0$), and let
\beqn{CRS-skdef}
s_k = \theta_k \frac{\epsilon^{1/2}}{2f_k}
\ms
(\theta_k > 0)
\eeqn
be an acceptable step for an arbitrary method in the CRS$_a$ class.
Now, because of \req{CRS-sigmakb}, \req{CRS-lambda-bounds} reduces to
\beqn{CRS-lambda2}
\lambda_k\in [0, \bar{\sigma}|s_k|]
= \left[0 , \bar{\sigma}\theta_k \frac{\epsilon^{1/2}}{2f_k}\right]
\eeqn
and, given that $H_k >0$ because of \req{CRS-fkgkHk}, this in turn implies
that $H_k+\lambda_k > 0$. Condition  \req{CRS-srule2} requires that
\beqn{CRS-rbound}
|g_k + (H_k+\lambda_k)s_k|
= |r_k|
\leq \kappa_{rg} |g_k|
= 2\kappa_{rg}\epsilon f_k
< 2 \epsilon,
\eeqn
where we used the fact that $f_k\leq 1$ because of \req{CRS-fkgkHk} and
$\kappa_{rg}< 1$ because of \req{CRS-srule2}. Moreover, \req{CRS-rbound}
and \req{CRS-lambda2} imply that
\beqn{CRS-sbounds}
\frac{2(1-\kappa_{rg})\epsilon f_k}{4\epsilon^{1/2}f_k^2 + \bar{\sigma}s_k}
\leq \frac{|g_k|(1-\kappa_{rg})}{H_k + \lambda_k}
\leq s_k
\leq \frac{|g_k|(1+\kappa_{rg})}{H_k + \lambda_k}
\leq \frac{(1+\kappa_{rg})\epsilon^{1/2}}{2 f_k}.
\eeqn
Thus, using \req{CRS-skdef} and the right-most part of these inequalities, we
obtain that $\theta_k \leq 1+\kappa_{rg}$, which in turn ensures that
$s_k \leq (1+\kappa_{rg})\epsilon^{1/2}/(2f_k)$.  Substituting this latter bound
in the denominator of the left-most part of \req{CRS-sbounds} and using
\req{CRS-skdef} again with the fact that $f_k \geq \half$ before termination, we
obtain that
\beqn{CRS-theta-bounds}
\theta_k
\in \left[ \frac{1-\kappa_{rg}}{1+\bar{\sigma}(1+\kappa_{rg})},1+\kappa_{rg}\right]
\eeqn
(note that this is \req{ex-skdef} with $\kappa_\lambda = 1+\bar{\sigma}(1+\kappa_{rg})$).
We immediately note that $\pi_k$ and $\phi(\theta_k)$ are then both guaranteed
to be bounded above and below as in \req{ex-piphi-def}. (Since this is enough
for our purpose, we ignore the additional restriction on $\theta_k$ which might
result from \req{CRS-slope}.)
\comment{%%%%%%%%%%%%%%%%%%%%%%%%%%%%%%%%%%%%%%%%%%%%%%%%%%%%%%%%%%%%
Now \req{CRS-rbound} gives that condition \req{CRS-slope},i.e.\ 
\beqn{CRS-slope2}
r_k
\leq \theta_k \frac{\epsilon^{1/2}}{4f_k}(
               4\epsilon^{1/2}f_k^2 +\lambda_k)
               + \frac{\kappa_1\theta_k^2}{8f_k^2}\epsilon
\leq \epsilon \theta_k \left[ \left(4 +\frac{\bar{\sigma}}{8f_k^2}\right)
               +\frac{\kappa_1\theta_k}{8f_k^2}\right],
\eeqn
is satisfied when $\theta_k \geq \half$. As a consequence, we deduce that the
step $s_k$ given by \req{CRS-skdef} is acceptable within the CRS$_a$ framework
provided
\beqn{CRS-theta-bounds}
\theta_k \in \left[ \,\max\left[\, \frac{1}{2},\,
      \frac{1}{1+\bar{\sigma}(1+\kappa_{rg})} \,\right],\,1+\kappa_{rg}\,\right].
\eeqn
We immediately note that, for such a step, $\pi_k$ and $\phi(\theta_k)$ are
both bounded above and below as in \req{ex-piphi-def}.
}%%%%%%%%%%%%%%%%%%%%%%%%%%%%%%%%%%%%%%%%%%%%%%%%%%%%%%%%%%%%%%%%%
Using the definitions
\req{CRS-fkgkHk} for $k+1$, we may then construct the objective function
$f_\epsilon^{CRS}$ on the interval $[x_k, x_k+s_k]$ by Hermite interpolation,
as in Section~\ref{slowMa-s}. Moreover, using \req{CRS-rho}, \req{CRS-fkgkHk},
\req{CRS-skdef}, \req{CRS-theta-bounds}, $f_k \in [\half,1]$ and the condition
\req{CRS-eta-cond}, we obtain that 
\[
\rho_k
=\frac{\epsilon^{3/2}}{2}\left(\frac{2f_k}{\theta_k \epsilon^{1/2}}\right)^3
=\frac{4 f_k^3}{\theta_k^3}
\geq \frac{1}{2(1+\kappa_{rg})^3}
\geq \eta.
\]
Thus iteration $k$ is successful, $x_{k+1}=x_k+s_k$, $\sigma^L_{k+1} = \sigma^L_k = 0$,
$\sigma^U_{k+1} = \sigma^U_k = \bar{\sigma}$, and all subsequent iterations
of the CRS$_a$ method up to termination follow the same pattern in accordance
with \req{CRS-fkgkHk}. As in Section~\ref{slowMa-s}, we may construct
$f_\epsilon^{CRS}$ on the whole of $\Re$ which satisfies A.1 and such that,
the considered CRS$_a$ method  applied to $f_\epsilon^{CRS}$  will terminate in
exactly $\lceil \epsilon^{-3/2} \rceil$ iterations and evaluations. This
and the $\calO\left(\epsilon^{-3/2}\right)$ upper bound on the worst-case
evaluation complexity of CRS methods allow stating the following theorem.

\lthm{slowCRS}{ For every $\epsilon \in (0,1)$ and
  every method in the CRS$_a$ class, a function $f^{CRS}_\epsilon$ satisfying
  A.$1$ with values in a bounded interval independent of $\epsilon$ can be
  constructed, such that the considered method terminates 
  exactly at iteration
  \[
  k_\epsilon=\left\lceil\epsilon^{-3/2}\right\rceil
  \]
  with the first iterate $x_{k_\epsilon}$ such that $\|\nabla_x f^{CRS}_\epsilon(x_{k_\epsilon})\| \leq
  \epsilon$. As a consequence, methods in CRS$_a$ are optimal within the CRS
  class and their worst-case evaluation complexity is, in order, also optimal
  with respect to that of methods in $\calM.1$.
}

CRS$_a$ then constitutes a kernel of optimal  methods (from the worst-case
evaluation complexity point of view) within CRS and $\calM.1$. Methods in CRS but
not in CRS$_a$ correspond to very inaccurate minimization of the regularized
model, which makes it unlikely that their worst-case evaluation complexity
surpasses that of methods in CRS$_a$. Finally note that, since we did not use
\req{CRS-slope} to construct our example, it effectively applies to a class
larger than CRS$_a$ where this condition is not imposed.

\numsection{The algorithm of Royer and Wright}

We finally consider the linesearch algorithm proposed in \cite[Algorithm~1]{RoyeWrig17},
which is reminiscent of the double linesearch algorithm of
\cite{GoulLuciRomaToin98} and \cite[Section~10.3.1]{ConnGoulToin00}.
From a given iterate $x_k$, this algorithm computes a search direction $d_k$
whose nature depends on the curvature of the (unregularized) quadratic model
along the negative gradient, and possibly computes the left-most eigenpair of
the Hessian if this curvature is negative or if the gradient's norm is small
enough to declare first-order stationarity. A linesearch along $d_k$ is then
performed by reducing the steplength $\alpha_k$ from $\alpha_k=1$ until
\beqn{RW-ls}
f(x_k+\alpha_k d_k) \leq f(x_k) - \frac{\eta}{6}\alpha_k^3 \|d_k\|^3
\eeqn
for some $\eta >0$.
The algorithm uses $\epsilon_g$ and $\epsilon_H$, two
different accuracy thresholds for first- and second-order approximate
criticality, respectively.

Our objective is now to show that, when applied to the function
$f^{\calM.1}_{\epsilon_g}$ of Section~\ref{slowMa-s} with $\epsilon = \epsilon_g$,
this algorithm, which we call the RW algorithm, takes exactly
$k_{\epsilon_g,1}= \lceil\epsilon_g^{-3/2}\rceil$ iterations and evaluations
to terminate with $\|g_k\| \leq \epsilon_g$. 

We first note that \req{ex-fkgkHk} guarantees that $H_k$ is positive definite
and, using \req{ex-fkbounds}, that
\[
\frac{g_k^TH_kg_k}{\|g_k\|^2} = 4 \epsilon_g^{1/2}f_k^2 > \epsilon_g
\]
for $k \in \iibe{0}{k_{\epsilon_g,1}}$. Then, provided
\beqn{RW-epsH-cond}
\epsilon_H \leq \sqrt{\epsilon_g},
\eeqn
and because $\lambda_{\min}(H_k) = 4 \epsilon_g^{1/2}f_k^2 > \epsilon_H$
(using \req{ex-fkbounds} again), the
RW algorithm defines the search direction from  Newton's equation $H_k d_k =
-g_k$ (which corresponds, as we have already seen, to taking $M_k=0=r_k$ and
thus $\theta_k=1$ in the example of Section~\ref{slowMa-s}). The RW algorithm
is therefore, on that example, identical to a linesearch variant of Newton's
method with the specific linesearch condition \req{RW-ls}.  Moreover, using
\req{ex-fkbounds} once more,
\[
f(x_k)-f(x_k+d_k)
= \frac{1}{2} \epsilon_g^{3/2}
\geq \frac{\eta}{6} \left( \frac{\epsilon_g^{1/2}}{2f_k}\right)^3
\geq \frac{\eta}{6} \epsilon_g^{3/2}
\]
whenever $\eta \leq 3$, an extremely weak condition\footnote{In practice,
  $\eta$ is most likely to belong to $(0,1)$ and even be reasonably close to
  zero.}.  Thus \req{RW-ls} holds\footnote{But fails for the example of
  Section~\ref{slowN-s} as $\|s_k\|=1$.} with $\alpha_k=1$. We have thus proved that
the RW algorithm generates the same sequence of iterates as Newton's method
when applied to $f^{\calM.1}_{\epsilon_g}$. The fact that an upper bound of
$\calO\left(\epsilon_g^{-3/2}\right)$ iterations and evaluations was proved to hold
in \cite[Theorem~5]{RoyeWrig17} then leads us to stating the following result.

\lthm{slowRW}{Assume that $\eta \in (0,3]$.  Then, for every $\epsilon_g \in
  (0,1)$ and $\epsilon_H$ satisfying \req{RW-epsH-cond}, a function
  $f^{\calM.1}_{\epsilon_g}$ satisfying A.$1$ with values in a bounded interval
  (independent of $\epsilon_g$ and $\epsilon_H$) can be constructed, such that
  the Royer-Wright algorithm terminates exactly at iteration
  \[
  k_{\epsilon_g}=\left\lceil\epsilon_g^{-3/2}\right\rceil
  \]
  with the first iterate $x_{k_{\epsilon_g}}$ such that
  $\|\nabla_x f^{\calM.1}_{\epsilon_g}(x_{k_{\epsilon_g}})\| \leq \epsilon_g$.
  As a consequence and under assumption \req{RW-epsH-cond}, the
  first-order worst-case evaluation complexity order of
  $\calO\left(\epsilon_g^{-3/2}\right)$ for this algorithm is sharp and it is
  (in order of $\epsilon_g$), also optimal with respect to that of algorithms
  in the $\calM.1$ and CRS classes.
}

\numsection{Conclusions}

We have provided lower bounds on the worst-case evaluation complexity of a
wide class of second-order methods for reaching approximate first-order
critical points of nonconvex, adequately smooth unconstrained optimization
problems. This has been achieved by providing improved examples of slow
convergence on functions with bounded range independent of $\epsilon$. We have
found that regularization algorithms, methods belonging to a subclass of that
proposed in \cite{CurtRobiSama17c} and the linesearch algorithm of
\cite{RoyeWrig17} are optimal from a worst-case complexity point of view
within a very wide class of second-order methods, in that their upper
complexity bounds match in order the lower bound we have shown for relevant,
sufficiently smooth objectives satisfying A.$\alpha$. At this point, the
question of whether all known optimal second-order methods share enough design
concepts to be made members of a single class remains open.

Note that every iteration complexity bound discussed above is of the order
$\epsilon^{-p}$ (for various values of $p>0$) for driving the objective's
gradient below $\epsilon$; thus the methods we have addressed may require an
exponential number of iterations $10^{p\cdot k}$ to generate $k$ correct
digits in the solution. Also, as our examples are one-dimensional, they fail
to capture the problem-dimension dependence of the upper complexity
bounds. Indeed, besides the accuracy tolerance $\epsilon$, existing upper
bounds depend on the distance to the solution set, that is $f(x_0)-f_{\rm
low}$, and the gradient's and Hessian's Lipschitz or H\"older constants, all
of which may dependent on the problem dimension. Some recent developments in
this respect can be found in~\cite{Jarr13,Agaretal16,JianLinMaZhan16,RoyeWrig17}.

Here we have solely addressed the evaluation complexity of generating
first-order critical points, but it is common to require second-order methods
for nonconvex problems to achieve second-order criticality.  Indeed, upper
worst-case complexity bounds are known in this case for cubic regularization
and trust-region methods \cite{NestPoly06, CartGoulToin11d, CartGoulToin12d},
which are essentially sharp in some cases \cite{CartGoulToin12d}.  A lower
bound on the whole class of second order methods for achieving second-order
optimality remains to be established, especially when different accuracy is
requested in the first- and second-order criticality conditions.

Regarding the worst-case evaluation complexity of constrained optimization
problems, we have shown \cite{CartGoulToin12b,CartGoulToin11a,CartGoulToin11b}
that the presence of constraints does not change the order of the bound, so
that the unconstrained upper bound for some first- or second-order methods
carries over to the constrained case; note that this does not include the cost
of solving the constrained subproblems as the latter does not require
additional problem evaluations.  Since constrained problems are at least as
difficult as unconstrained ones, these bounds are also sharp. It remains an
open question whether a unified treatment such as the one given here can be
provided for the worst-case evaluation complexity of methods for constrained
problems.

\appendix

\appnumsection{A1. An example of slow convergence of the steepest-descent method}

We show in this paragraph that the steepest-descent method may need at least
$\epsilon^{-2}$ iteration to terminate on a function whose range is fixed and
independent of $\epsilon$.

We once again follow the methodology used in Section~\ref{slowMa-s} and
build a unidimensional function $f^{SD}_\epsilon$ by Hermite interpolation,
such that the steepest-descent method applied to this function takes exactly
$k_{\epsilon} = \lceil\epsilon^{-2}\rceil$
iterations and function evaluations to terminate with an iterate $x_k$ such
that $|g(x_k)| \leq \epsilon$.  Note that, for the sequence of function values
to be interpretable as the result of applying the steepest-descent method (using a
Goldstein linesearch), we require that, for all $k$,
\beqn{is-LS}
 f(x_k) + \mu_1 g_k^Ts_k \leq  f(x_k-\mu_kg_k) \leq f(x_k) + \mu_2 g_k^Ts_k
\tim{ for constants }
0 < \mu_2 < \mu_1 < 1
\eeqn
where, as above, $s_k = x_{k+1}-x_k$. Keeping this in mind, we define the
sequences ${f_k}$, ${g_k}$, ${H_k}$ and ${s_k}$ for $k \in
\iibe{0}{k_{\epsilon}-1}$ by
\[
f_k = 1 - \half k \epsilon^2
\ms
g_k = - 2 \epsilon f_k,
\ms
H_k =  0,
\ms
r_k = 0
\tim{and}
\mu_k = \frac{1}{4 f_k^2} \in [\quarter, 1].
\]
Note that this last definition ensures that \req{is-LS} holds provided $0
<\mu_2 < \half < \mu_1 < 1$. It also gives that
$s_k = \epsilon/(2 f_k) \leq \epsilon < 1$.
Using these values, it can also be verified that
termination occurs for $k = k_{\epsilon}$, that
$f^{SD}_\epsilon$ defined by \req{fM1} and Hermite interpolation is twice
continuously differentiable on $[0,x_{k_\epsilon}]$ and that \req{ex-c0c1}
again holds. Since $|g_k| \leq \epsilon$, we also obtain that, for $k \in
\iibe{0}{k_{\epsilon}-1}$,
\[
\left|\frac{\Delta f_k}{s_k^2}\right|
= 2 f_k^2 \leq 1,
\ms
\left|\frac{\Delta g_k}{s_k}\right|
= 2 \epsilon^2 f_k \leq 2
\tim{and}
\left|\frac{g_k}{s_k}\right|
= 4 f_k^2 \leq 4.
\]
These bounds, $H_k = \Delta H_k= 0$, the first equality of
\req{bounded-Hessian} and \req{ex-expck} then imply that the Hessian of
$f^{SD}_\epsilon$ is bounded above by a constant independent of $\epsilon$.
$f^{SD}_\epsilon$ thus satisfies A.$0$ and therefore has Lipchitz continuous
gradient.  Moreover, since $s_k \leq 1$, we also obtain, as in
Section~\ref{slowMa-s} and \ref{slowN-s}, that $|f^{SD}_\epsilon|$ is bounded
by a constant independent of $\epsilon$ on $[0,x_{k_\epsilon}]$.  As above we 
then extend $f^{SD}_\epsilon$ to the whole of $\Re$ while preserving A.$0$.

\lthm{slowSD}{ For every $\epsilon \in (0,1)$, a function
  $f^{SD}_\epsilon$ satisfying A.$0$ (and thus having Lipschitz continuous
  gradient) with values in a bounded interval independent of $\epsilon$ can be
  constructed, such that the steepest-descent method terminates exactly at
  iteration
  \[
  k_\epsilon=\left\lceil\epsilon^{-2}\right\rceil
  \]
  with the first iterate $x_{k_\epsilon}$ such that
  $\|\nabla_xf^{SD}_\epsilon(x_{k_\epsilon})|\leq \epsilon$.
}

As a consequence, the $\calO\left(\epsilon^{-2}\right)$ order of worst-case evaluation
complexity is sharp for the steepest-descent method in the sense that the
complexity ratio $\rho_{\rm comp}$ is bounded above independently of of
$\epsilon$, which improves on the conclusion proposed in
\cite{CartGoulToin10a} for the steepest-descent method.

\begin{figure}[htbp] % produced by exsd.m
\begin{center}
\vspace*{1.5mm}
\includegraphics[height=3.5cm]{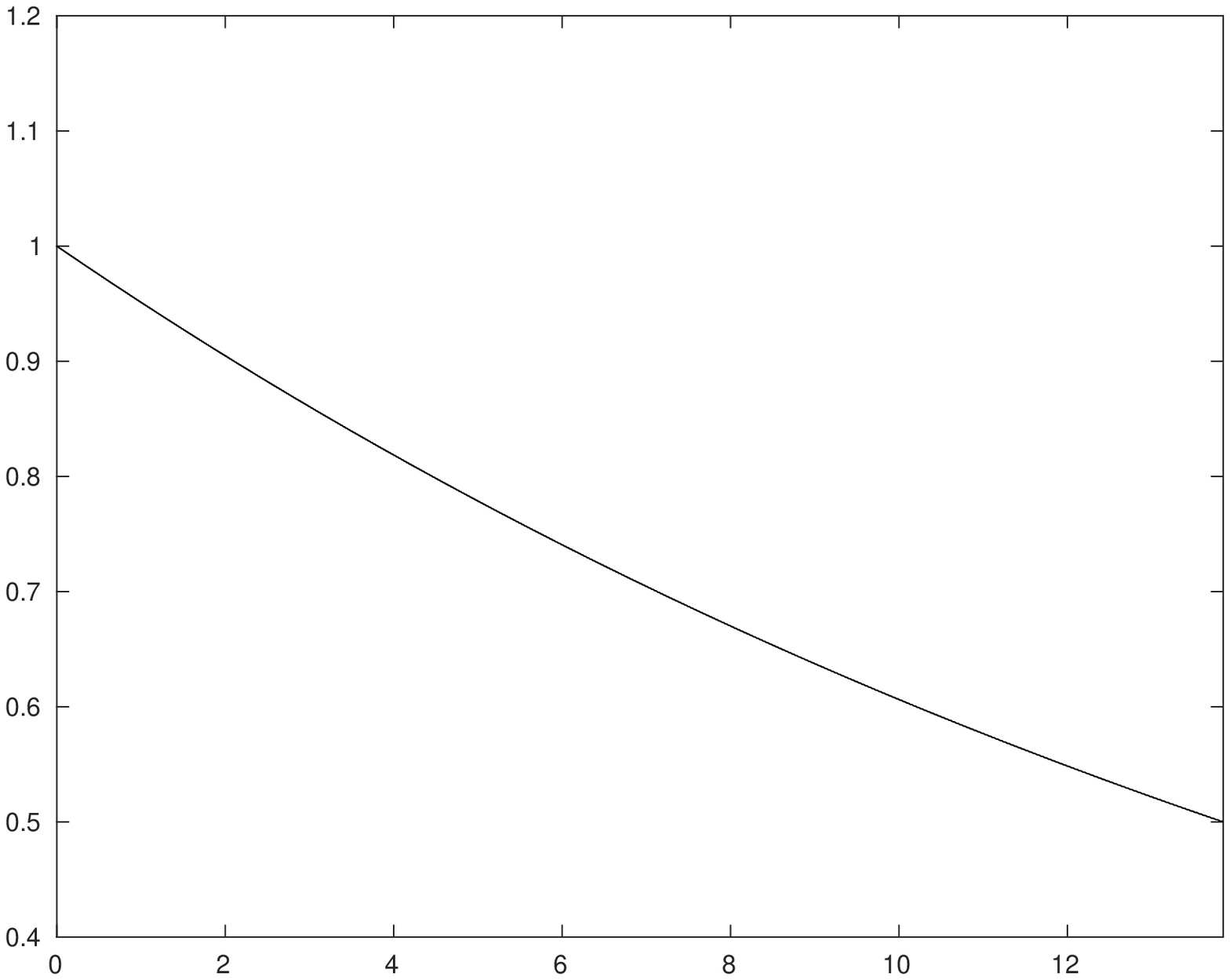}  
\hspace*{4mm}
\includegraphics[height=3.5cm]{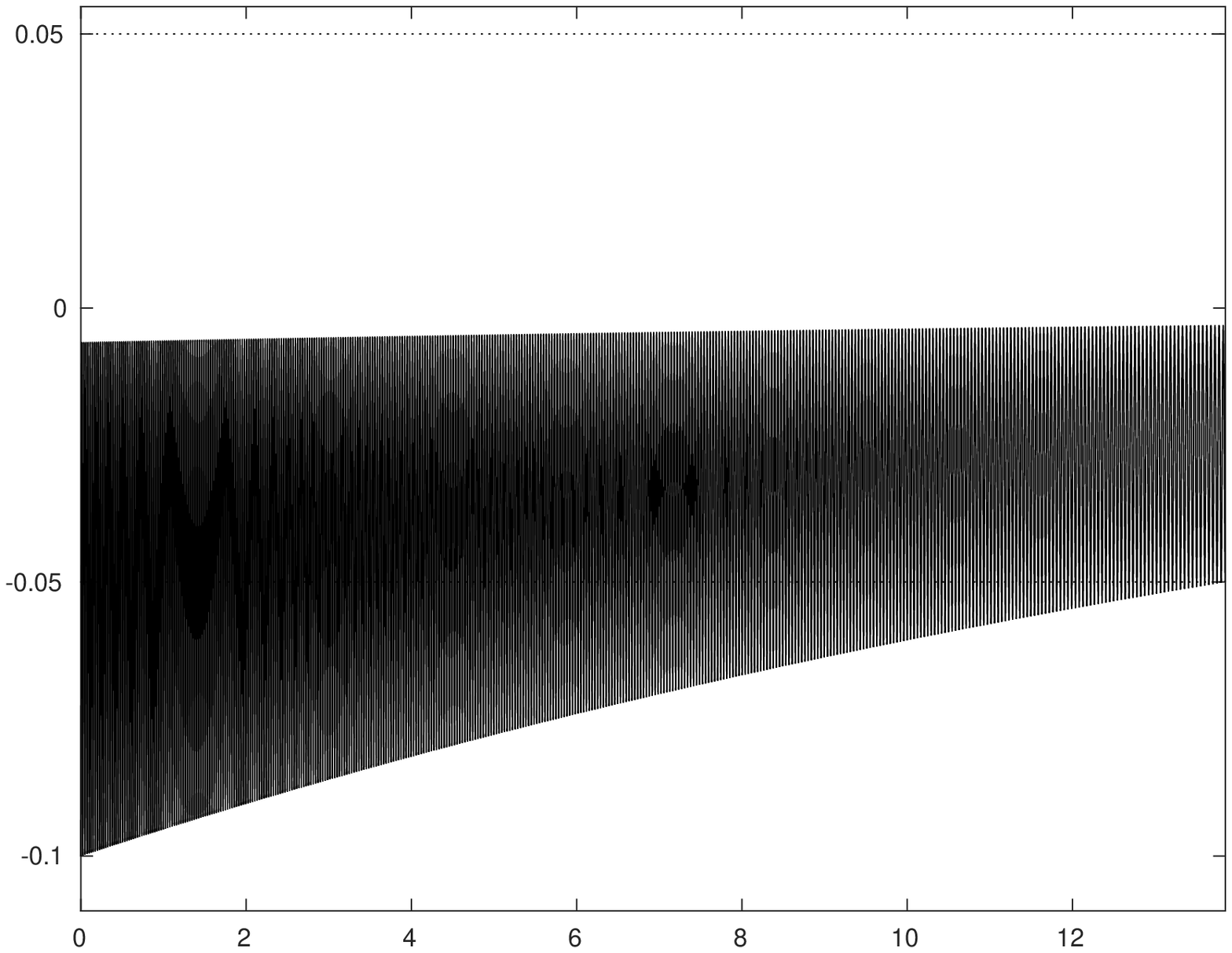}
\hspace*{4mm}
\includegraphics[height=3.5cm]{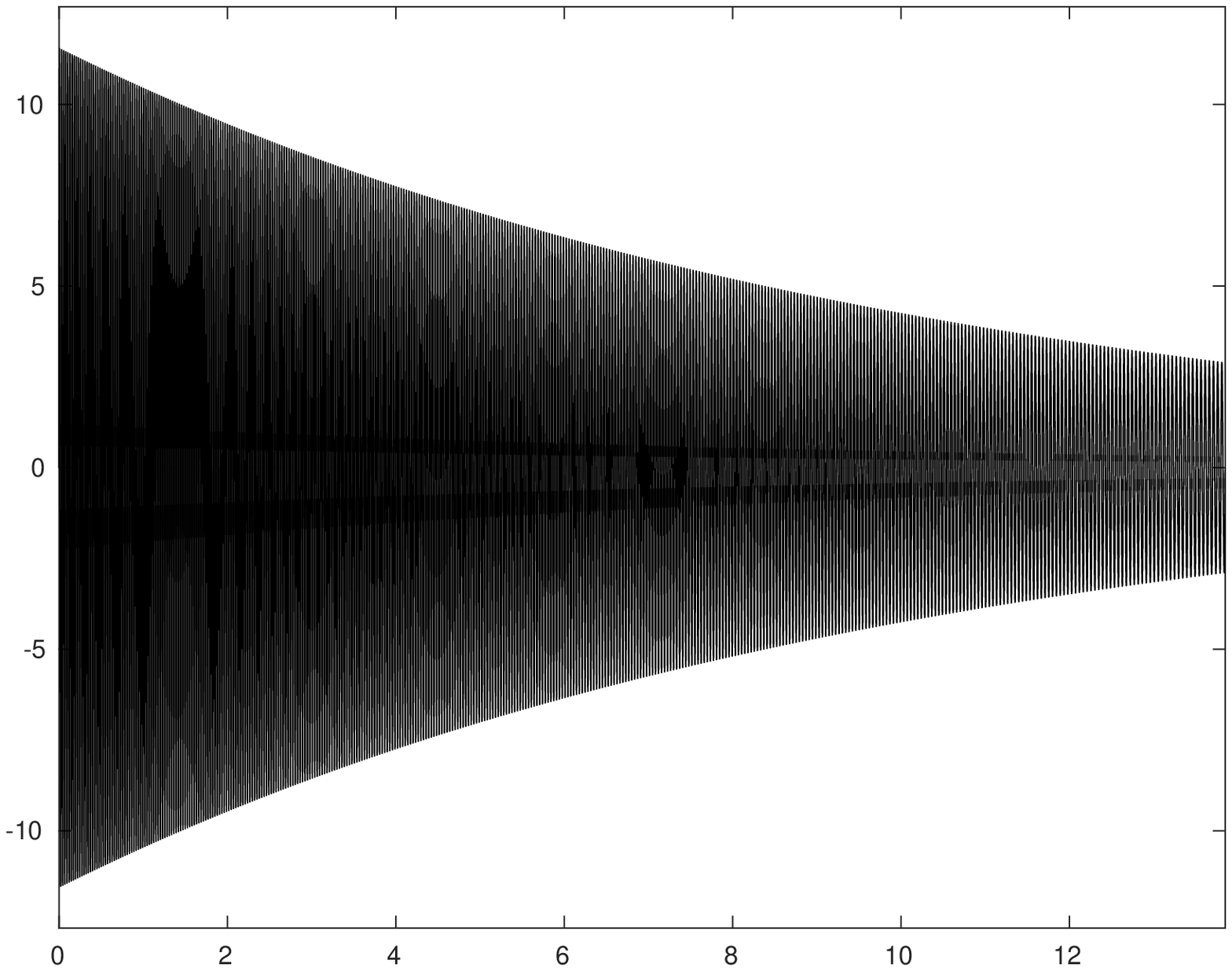}\\*[1.5ex]
\includegraphics[height=3.5cm]{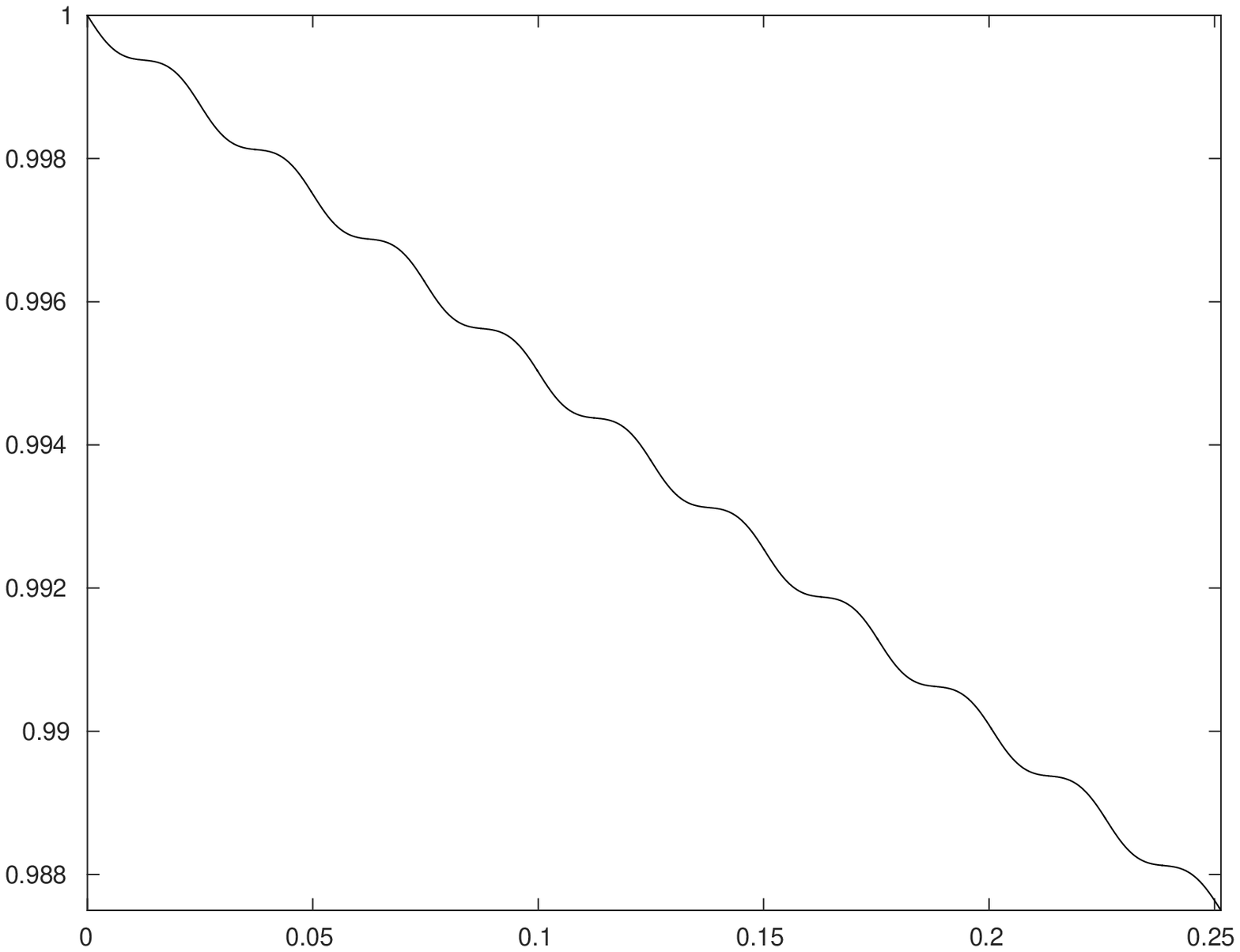}  
\hspace*{4mm}
\includegraphics[height=3.5cm]{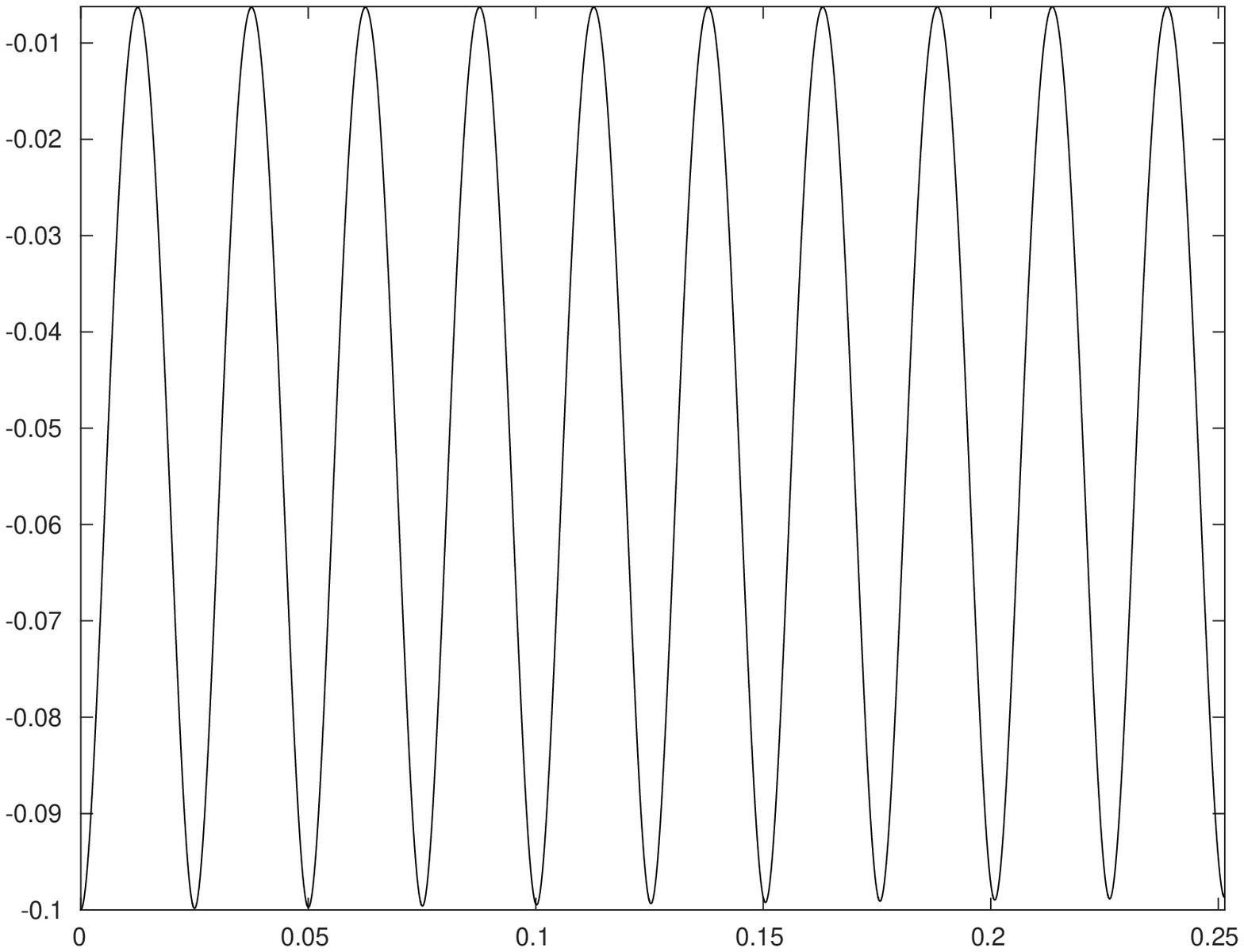}
\hspace*{4mm}
\includegraphics[height=3.5cm]{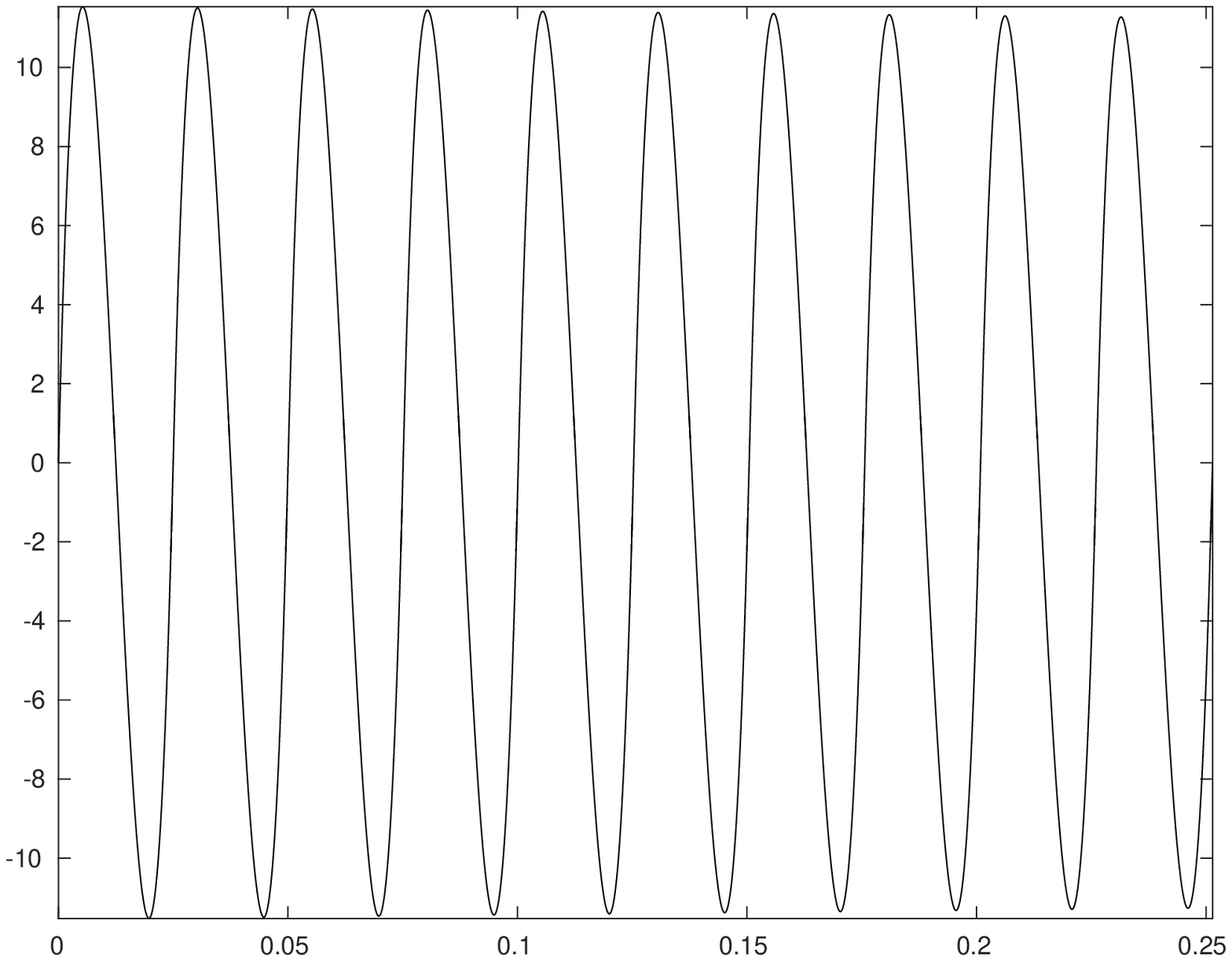}
\end{center}
\vspace*{-5mm}
\caption{\label{exsd}$f^{SD}_\epsilon(x)$ (left) and its first (center) and
  second (right) derivatives as a function of $x$ for $\epsilon = 5.10^{-2}$
  (top: $x \in [0,x_{k_{\epsilon,\alpha}}]$; bottom: $x \in [0,x_{10}]$).
  Horizontal dotted lines indicate values of $-\epsilon$ and $\epsilon$ in the
  central top graph.}
\end{figure}

The top three graphs of Figure~\ref{exsd} illustrate the global behaviour
of the resulting function $f^N_\epsilon(x)$ and of its first and second
derivatives for $x \in [0,x_{k_\epsilon}]$, while the bottom ones show more detail of the first
10 iterations. The figure is once more constructed using $\epsilon =
5.10^{-2}$ ($k_{\epsilon} = 400$).

\appnumsection{A2. Upper complexity bound for the $(2+\alpha)$-regularization method}

The purpose of this paragraph is to to provide some of the missing details in
the proof of Lemma~\ref{lemma:regularization}, as well as making explicit the
statement made at the end of Section~5.1 in \cite{CartGoulToin11d} that the
$(2+\alpha)$-regularization method needs at most \req{upper-reg} iterations
(and function/derivatives evaluations) to obtain and iterate $x_k$ such that
$|g_k| \leq \epsilon$.

We start by proving \req{reg-sk-upper} following the
reasoning of \cite[Lem.2.2]{CartGoulToin11}. Consider
\[
\begin{array}{lcl}
m_k(s) - f(x_k)
& = & g_k^Ts + \half s^TH_ks + \frac{1}{2+\alpha}\sigma_k \|s\|^{2+\alpha} \\*[1ex]
& \geq & -\|g_k\|\,\|s\| - \half \|s\|^2\, \|H_k\| + \frac{1}{2+\alpha}\sigma_k \|s\|^{2+\alpha} \\*[1ex]
& \geq & \left(\frac{1}{3(2+\alpha)} \sigma_k \|s\|^{2+\alpha}-\|g_k\|\,\|s\|\right)
         + \left( \frac{2}{3(2+\alpha)} \sigma_k \|s\|^{2+\alpha}- \half \|s\|^2\|H_k\|\right)
\end{array}
\]
But then $\frac{2}{3(2+\alpha)} \sigma_k \|s\|^{2+\alpha}-\|H_k\|\,\|s\|^2 >0$
if
%\[
%\|s_k\| < \left(\frac{3(2+\alpha)\|H_k\|}{4\sigma_k}\right)^{\frac{1}{\alpha}}
%\]
$\|s_k\| < (3(2+\alpha)\|H_k\|/(4\sigma_k))^{\frac{1}{\alpha}}$
while $\frac{1}{3(2+\alpha)} \sigma_k \|s\|^{2+\alpha}-\|g_k\|\,\|s\|> 0$ if
%\[
%\|s_k\| < \left(\frac{3(2+\alpha)\|g_k\|}{\sigma_k}\right)^{\frac{1}{1+\alpha}}.
%\]
$\|s_k\| < (3(2+\alpha)\|g_k\|/\sigma_k)^{\frac{1}{1+\alpha}}$.
Hence, since $m_k(s_k) < f(x_k)$, we have that
\[
\|s_k\| \leq
\max\left[\left(\frac{3(2+\alpha)\|H_k\|}{4\sigma_k}\right)^{\frac{1}{\alpha}},
   \left(\frac{3(2+\alpha)\|g_k\|}{\sigma_k}\right)^{\frac{1}{1+\alpha}}
   \right]
\]
which yields \req{reg-sk-upper} because $\|H_k\| \leq L_g$.

We next explicit the worst-case evaluation complexity bounf of Section~5.1 in
\cite{CartGoulToin11d}.  Following
\cite[Lemma~5.2]{CartGoulToin11}, we start by proving that
\beqn{sigma-max}
\sigma_{\max}\eqdef c_{\sigma}\max(\sigma_0,L_{H,\alpha})
\eeqn
for some constant $c_\sigma$ only dependent on $\alpha$ and algorithm's
parameters. To show this inequality, we deduce from Taylor's theorem that, for
each $k\geq 0$ and some $\xi_k$ belonging the the segment $[x_k, x_k+s_k]$,
\[
f(x_k+s_k)-m_k(s_k)
\leq \frac{1}{2}\|H(\xi_k)-H(x_k)\|\cdot\|s_k\|^2-\frac{\sigma_k}{2+\alpha}\|s_k\|^{2+\alpha}
\leq \left(\frac{L_{H,\alpha}}{2}-\frac{\sigma_k}{2+\alpha}\right)\|s_k\|^{2+\alpha},
\]
where, to obtain the second inequality, we employed \req{LipsHNEW} in 
A.$\alpha$ and $\|\xi_k-x_k\|\leq \|s_k\|$. Thus $f(x_k+s_k) < m_k(s_k)$
whenever $\sigma_k > \half(2+\alpha)L_{H,\alpha}$, providing sufficient
descent and ensuring that $\sigma_{k+1} \leq \sigma_k$.  Taking into account
the (possibly large) choice of the regularization parameter at startup then
yields \req{sigma-max}. 

We next note that, because of \req{lbCAL} and \req{sigma-max}, \req{Msk}
holds. Moreover, $\kappa(M_k)= \kappa\left(\sigma_k\|s_k\|^\alpha I\right) =
1$. Lemma~\ref{steppropertyLemma} then ensures that \req{sklbd} also holds. 

We finally follow \cite[Corollary 5.3]{CartGoulToin11} to prove the final upper
bound on the number of successful iterations (and hence on the number of
function and derivatives evaluations). Let $\calS^\epsilon_k$ index
the subset of the first $k$ iterations that are successful and such that
$\min[\|g_k\|,\|g_{k+1}\|] > \epsilon$, and let $|\calS^\epsilon_k|$
denote its cardinality. It  follows  from this definition, \req{Msk},
\req{sigmamin} and the fact that sufficient decrease is obtained at
successful iterations that, for all $k$ before termination,
\beqn{finalbdmk}
f(x_j)-m_k(s_j)\geq \alpha_{\rm S}\epsilon^{\frac{2+\alpha}{1+\alpha}},
\tim{for all $j\in \calS^\epsilon_k$,}
\eeqn
for some positive constant $\alpha_{\rm S}$ independent of $\epsilon$.
Now, if $f_{\rm low} > -\infty$ is a lower bound on $f(x)$, we have, using the
monotonically decreasing nature of $\{f(x_k)\}$, that
\[
\begin{array}{lcccl}
f(x_0)-f_{\rm low}
& \geq &f(x_0)-f(x_{k+1}) %\\*[1ex]
& = & \bigsum_{j \in \calS^\epsilon_k} \left[f(x_j)-f(x_{j+1})\right]\\*[1ex]
& \geq & \eta_1 \bigsum_{j \in \calS^\epsilon_k}\left[f(x_j)-m_k(s_j)\right]%\\*[1ex]
& \geq & |\calS^\epsilon_k| \,\eta_1 \alpha_{\rm S}\,\epsilon^{\frac{2+\alpha}{1+\alpha}},
\end{array}
\]
where the constant $\eta_1 \in (0,1)$ defines sufficient decrease.
Hence, for all $k \geq 0$,
\[
|\calS^\epsilon_k|
\leq \frac{f(x_0)-f_{\rm low}}{\eta_1 \alpha_{\rm S}} \,\, \epsilon^{-\frac{2+\alpha}{1+\alpha}}.
\]
As a consequence, the $(2+\alpha)$-regularization method needs at most
\req{upper-reg} successful iterations to terminate.  Since it known that,
for regularization methods, $k \leq \kappa_\calS |\calS^\epsilon_k|$
for some constant $\kappa_\calS$ \cite[Theorem~2.1]{CartGoulToin11d} and
because every iteration involves a single evaluation, we conclude that the
$(2+\alpha)$-regularization method needs at most \req{upper-reg} function and
derivatives evaluations to produce an iterate $x_k$ such that $\|g_k\| \leq
\epsilon$ when applied to an objective function satisfying A.$\alpha$.

We finally oserve that the statement (made in the proof of
Lemma~\ref{lemma:regularization}) that $\|g_k\|$ is bounded above immediately
follows from this worst-case evaluation complexity bound.

\end{document}